\documentclass[10pt]{article}

\usepackage{amsmath, amsthm, amssymb, enumerate}
\usepackage{color}
\usepackage{CJK, times}
\usepackage{datetime}
\usepackage{fancyhdr}

\usepackage[backgroundcolor=gray!30,linecolor=black]{todonotes}
\chardef\coloryes=1 %%%out
\chardef\isitdraft=0 %%%out
\ifnum\isitdraft=1 %%%out
   \input macros0.tex %%%out
   \def\eqref#1{({\ref{#1}})}                %saves writing paranthesis%%%out

\usepackage[color,notcite]{showkeys}%%%out
\definecolor{refkey}{gray}{.3}%%%out

\definecolor{labelkey}{gray}{.3}%%%out

\else%%%out

\fi%%%out

\newcommand{\ip}[2]{\left<#1,#2\right>}
\begin{document}
%\newcommand{\llabel}{\label}             %\llabel is just a synonim for \label
%\newcommand{\rref}{\ref}                 %\rref is  just a synonim for \ref
%\newcommand{\ccite}{\cite}               %\ccite is just a synonim for
                                         %when ref. to equations
\def\Dg{{D'g}}
\def\ua{u^{\alpha}}
\def\intint{\int\!\!\!\!\int}
\def\intinttext{\int\!\!\!\int}
\def\intintint{\int\!\!\!\!\int\!\!\!\!\int}
\def\intintintint{\int\!\!\!\!\int\!\!\!\!\int\!\!\!\!\int}
\def\ques{{\cor \underline{??????}\cob}}
\def\nto#1{{\coC \footnote{\em \coC #1}}}
\def\fractext#1#2{{#1}/{#2}}
\def\fracsm#1#2{{\textstyle{\frac{#1}{#2}}}}   %smaller version of frac
\def\baru{U}
\def\nnonumber{}
\def\palpha{p_{\alpha}}
\def\valpha{v_{\alpha}}
\def\qalpha{q_{\alpha}}
\def\walpha{w_{\alpha}}
\def\falpha{f_{\alpha}}
\def\dalpha{d_{\alpha}}
\def\galpha{g_{\alpha}}
\def\halpha{h_{\alpha}}
\def\psialpha{\psi_{\alpha}}
\def\psibeta{\psi_{\beta}}
\def\betaalpha{\beta_{\alpha}}
\def\gammaalpha{\gamma_{\alpha}}
\def\Talpha{T}
\def\TTalpha{T_{\alpha}}
\def\TTalphak{T_{\alpha,k}}
\def\falphak{f^{k}_{\alpha}}

\def\cof{\mathop{\rm cf\,}\nolimits}
\def\abar{\overline a}
\def\vbar{\overline v}
\def\qbar{\overline q}
\def\wbar{\overline w}
\def\etabar{\overline\eta}
\def\Om{\Omega}
\def\ep{\epsilon}
\def\al{\alpha}
\def\Ga{\Gamma}
\def\R{\mathbb R}

\newcommand{\bv}{{u}}
\newcommand {\Dn}[1]{\frac{\partial #1  }{\partial \nu}}
\def\mm{m}

%if \isitdraft=0 or \coloryes=0%%%out
\def\cor{{}}%%%out
\def\cog{{}}%%%out
\def\cob{{}}%%%out
\def\coe{{}}%%%out
\def\coA{{}}%%%out
\def\coB{{}}%%%out
\def\coC{{}}%%%out
\def\coD{{}}%%%out
\def\coE{{}}%%%out
\def\coF{{}}%%%out
%\def\MR#1{}%%%out

%%3
%\ifnum\isitdraft=1%%%out
\ifnum\coloryes=1%%%out

  \definecolor{coloraaaa}{rgb}{0.1,0.2,0.8}%%%out
  \definecolor{colorbbbb}{rgb}{0.1,0.7,0.1}%%%out
  \definecolor{colorcccc}{rgb}{0.8,0.3,0.9}%%%out
  \definecolor{colordddd}{rgb}{0.0,.5,0.0}%%%out
  \definecolor{coloreeee}{rgb}{0.8,0.3,0.9}%%%out
  \definecolor{colorffff}{rgb}{0.8,0.3,0.9}%%%out
  \definecolor{colorgggg}{rgb}{0.5,0.0,0.4}%%%out

 \def\cog{\color{colordddd}}%%%out
 \def\cob{\color{black}}%%%out
 \def\cor{\color{red}}%%%out
 \def\coe{\color{colorgggg}}%%%out

 \def\coA{\color{coloraaaa}}%%%out
 \def\coB{\color{colorbbbb}}%%%out
 \def\coC{\color{colorcccc}}%%%out
 \def\coD{\color{colordddd}}%%%out
 \def\coE{\color{coloreeee}}%%%out
 \def\coF{\color{colorffff}}%%%out
 \def\coG{\color{colorgggg}}%%%out

%%4
\fi%%%out
 \ifnum\isitdraft=1%%%out
   \baselineskip=17pt%%%out
   \input macros.tex%%%out
   \def\blackdot{{\color{red}{\hskip-.0truecm\rule[-1mm]{4mm}{4mm}\hskip.2truecm}}\hskip-.3truecm}%%%out
   \def\bdot{{\coC {\hskip-.0truecm\rule[-1mm]{4mm}{4mm}\hskip.2truecm}}\hskip-.3truecm}%%%out
   \def\purpledot{{\coA{\rule[0mm]{4mm}{4mm}}\cob}}%%%out
   \def\purpledottwo{{\coA{\rule[0mm]{4mm}{4mm}}\cob}}%%%out
   \def\pdot{\purpledot}%%%out
\else%%%out  
   \baselineskip=15pt
   \def\blackdot{{\rule[-3mm]{8mm}{8mm}}}%%%out
   \def\purpledot{{\rule[-3mm]{8mm}{8mm}}}%%%out
   \def\pdot{}%%%out
%%5
\fi%%%out

\def\nts#1{{\hbox{\bf ~#1~}}} %nts=note to self
\def\nts#1{{\cor\hbox{\bf ~#1~}}} %nts=note to self%%%out
\def\ntsf#1{\footnote{\hbox{\bf ~#1~}}} %nts=note to self
\def\ntsf#1{\footnote{\cor\hbox{\bf ~#1~}}} %nts=note to self%%%out
\def\bigline#1{~\\\hskip2truecm~~~~{#1}{#1}{#1}{#1}{#1}{#1}{#1}{#1}{#1}{#1}{#1}{#1}{#1}{#1}{#1}{#1}{#1}{#1}{#1}{#1}{#1}\\}%%%out
\def\biglineb{\bigline{$\downarrow\,$ $\downarrow\,$}}%%%out
\def\biglinem{\bigline{---}}%%%out
\def\biglinee{\bigline{$\uparrow\,$ $\uparrow\,$}}%%%out

\def\mbar{{\overline M}}
\def\tilde{\widetilde}
\newtheorem{Theorem}{Theorem}[section]
\newtheorem{Corollary}[Theorem]{Corollary}
\newtheorem{Definition}[Theorem]{Definition}
\newtheorem{Proposition}[Theorem]{Proposition}
\newtheorem{Lemma}[Theorem]{Lemma}
\newtheorem{Remark}[Theorem]{Remark}
\newtheorem{definition}{Definition}[section]
\def\theequation{\thesection.\arabic{equation}}
\def\endproof{\hfill$\Box$\\}
\def\comma{ {\rm ,\qquad{}} }            %comma in a formula
\def\commaone{ {\rm ,\qquad{}} }         %second comma in a formula
\def\dist{\mathop{\rm dist}\nolimits}    %distance
\def\sgn{\mathop{\rm sgn\,}\nolimits}    %sgn
\def\Tr{\mathop{\rm Tr}\nolimits}    %trace
\def\div{\mathop{\rm div}\nolimits}    %divergence
\def\supp{\mathop{\rm supp}\nolimits}    %divergence
\def\divtwo{\mathop{{\rm div}_2\,}\nolimits}    %two dimensional divergence
\def\re{\mathop{\rm {\mathbb R}e}\nolimits}    %distance
\def\indeq{\qquad{}\!\!\!\!}                     %indentation in formulas
\def\period{.}                           %period in a formula
\def\semicolon{\,;}                      %semicolon in a formula
\newcommand{\cD}{\mathcal{D}}
%**end of header

\title{On the local well-posedness and a Prodi-Serrin type regularity criterion of the three-dimensional MHD-Boussinesq system without thermal diffusion}
\author{Adam Larios, Yuan Pei} 
\maketitle
\date{}
\bigskip
\indent Department of Mathematics\\
\indent University of Nebraska-Lincoln\\
\indent Lincoln, NE 68588\\
\indent e-mail: alarios\char'100unl.edu
\indent e-mail: ypei4\char'100unl.edu

\bigskip

\begin{abstract}

We prove a Prodi-Serrin-type global regularity condition for the three-dimensional Magnetohydrodynamic-Boussinesq system (3D MHD-Boussinesq) without thermal diffusion, in terms of only two velocity and two magnetic components. This is the first Prodi-Serrin-type criterion for a hydrodynamic system which is not fully dissipative, and indicates that such an approach may be successful on other systems. In addition, we provide a constructive proof of the local well-posedness of solutions to the fully dissipative 3D MHD-Boussinesq system, and also the fully inviscid, irresistive, non-diffusive MHD-Boussinesq equations. We note that, as a special case, these results include the 3D non-diffusive Boussinesq system and the 3D MHD equations. Moreover, they can be extended without difficulty to include the case of a Coriolis rotational term.

\end{abstract}

\noindent\thanks{\em Keywords:\/}
Magnetohydrodynamic equations, Boussinesq equations, 
B\'enard convection, Prodi-Serrin, 
weak solutions, partial viscosity, 
inviscid, global existence, regularity

\noindent\thanks{\em Mathematics Subject Classification\/}:
35A01, % Existence problems: global existence, local existence, non-existence
35K51, % Initial--boundary value problems for second-order parabolic systems
35Q35, % PDEs in connection with fluid mechanics
35Q86, % PDEs in connection with geophysics
76B03, % Existence, uniqueness, and regularity theory (Incompressible inviscid fluids)
76D03, % Existence, uniqueness, and regularity theory (Incompressible viscous fluids)
76W05.% Magnetohydrodynamics and electrohydrodynamics
\section{Introduction}
\label{sec1}
In this paper, we address global regularity criteria for the solutions to the
non-diffusive three-dimensional MHD-Boussinesq system of equations. 
The MHD-Boussinesq system models the convection of an incompressible flow 
driven by the buoyant effect of a thermal or density field, and the Lorenz force, generated by the magnetic field of the fluid. 
Specifically, it closely relates to a natural type of the Rayleigh-B\'enard convection, 
which occurs in a horizontal layer of conductive fluid heated from below, 
with the presence of a magnetic field 
(c.f. \cite{MP, MR}).
Various physical theories and numerical experiments such as in \cite{TM} have been developed 
to study the Rayleigh-B\'enard as well as the magnetic Rayleigh-B\'enard convection and related equations. 
We observe that by formally setting the magnetic field $b$ to zero, 
system (\ref{S2}) below reduces to the Boussinesq equations 
while by formally setting the thermal fluctuation $\theta=0$ we obtain the magnetohydrodynamic equations. 
One also formally recovers the incompressible Navier-Stokes equations if we set $b=0$ and $\theta=0$ simultaneously. 

Denote by $\Omega = \mathbb{T}^3$ the three-dimensional periodic space $\mathbb{R}^3/\mathbb{Z}^3 = [0, 1]^3$, 
and for $T>0$, the 3D MHD-Boussinesq system with full fluid viscosity, magnetic resistivity, and thermal diffusion 
over $\Omega\times [0, T)$ is given by
\begin{equation}
       \left\{
       \begin{aligned}
       &\frac{\partial u}{\partial t} 
       - \nu\Delta u 
       + (u\cdot\nabla) u 
       + \nabla p 
       = 
       (b\cdot\nabla) b
       +
       g\theta e_{3},
     \\&
       \frac{\partial b}{\partial t} 
       - \eta\Delta b 
       + (u\cdot\nabla) b
       = 
       (b\cdot\nabla) u,
     \\&
       \frac{\partial \theta}{\partial t} 
       -
        \kappa\Delta\theta
       + 
       (u\cdot\nabla) \theta  
       = 
       0, 
    \\&
    \nabla \cdot u = 0 = \nabla \cdot b,
    \end{aligned}
    \right.
    \label{S2}
  \end{equation}
where $\nu\geq0$, $\eta\geq0$, and $\kappa\geq0$ stand for the constant kinematic viscosity, magnetic diffusivity, and thermal diffusivity, respectively. The constant $g>0$ has unit of force, and is proportional to the constant of gravitational acceleration.  
We denote $x = (x_1, x_2, x_3)$, and $e_3$ to be the unit vector in the $x_3$ direction, i.e., $e_3 = (0, 0, 1)^{T}$.
Here and henceforth, $u = u(x, t) = (u_1(x, t), u_2(x,t), u_3(x, t))$ is the unknown velocity field of a viscous incompressible fluid, 
with divergence-free initial data $u(x, 0)=u_0$;
$b = b(x, t) = (b_1(x, t), b_2(x,t), b_3(x, t))$ is the unknown magnetic field, 
with divergence-free initial data $b(x, 0)=b_0 $;
and the scalar $p=p(x, t)$ represents the unknown pressure, while $\theta = \theta(x, t)$ can be thought of as the unknown temperature fluctuation, 
with initial value $\theta_{0} = \theta(x, 0)$. 
Setting $\kappa=0$, we obtain the non-diffusive MHD-Boussinesq system
\begin{equation}
       \left\{
       \begin{aligned}
       &\frac{\partial u}{\partial t} 
       - \nu\Delta u 
       +(u\cdot\nabla) u 
       + \nabla p 
       = 
       (b\cdot\nabla) b
       +
       g\theta e_{3},
     \\&
       \frac{\partial b}{\partial t} 
       - \eta\Delta b 
       + (u\cdot\nabla) b
       = 
       (b\cdot\nabla) u,
     \\&
       \frac{\partial \theta}{\partial t} 
       + 
       (u\cdot\nabla) \theta  
       = 
       0, 
    \\&
    \nabla \cdot u = 0 = \nabla \cdot b, 
    \end{aligned}
    \right.
    \label{S1}
  \end{equation}
which we study extensively in this paper. We also provide a proof for the local existence and uniqueness of solutions 
to the fully inviscid MHD-Boussinesq system with $\nu=\eta=\kappa=0$, namely, 
\begin{equation}
  \left\{
    \begin{aligned}
     &\frac{\partial u}{\partial t} 
       + (u\cdot\nabla) u 
       + \nabla p 
       = 
       (b\cdot\nabla) b
       +
       g\theta e_{3}, 
     \\&
       \frac{\partial b}{\partial t} 
       + (u\cdot\nabla) b
       = 
       (b\cdot\nabla) u, 
     \\&
       \frac{\partial \theta}{\partial t} 
       + 
       (u\cdot\nabla) \theta  
       = 
       0,
    \\&
    \nabla \cdot u = 0 = \nabla \cdot b,
   \end{aligned}
 \right.
 \label{S3}
  \end{equation}
with the initial condition $u_0$, $b_0$, and $\theta_0$ in $H^3$. We note that the proof of this result differs sharply from the proof of local existence for solutions of \eqref{S2}, due to a lack of compactness.  Therefore, we include the proof for the sake of completeness.

In recent years, from the perspective of mathematical fluid dynamics, 
much progress have been made in the study of solutions of the Boussinesq and MHD equations.
For instance, in \cite{ChaeKimNam1999, CN}, Chae et al. obtained the local well-posedness 
of the fully inviscid 2D Boussinesq equations with smooth initial data. 
A major breakthrough came in \cite{Ch} and \cite{HL}, where the authors independently proved global well-posedness for the two-dimensional Boussinesq equations with the case $\nu > 0$ and $\kappa=0$ and the case $\nu = 0$ and $\kappa > 0$, 
On the other hand, Wu et al. proved in \cite{CRW, CW1, CW2, DL, W2} 
the global well-posedness of the MHD equations, 
for a variety of combinations of dissipation and diffusion in two dimensional space.
Furthermore, a series of results concerning the global regularity 
of the 2D Boussinesq equations with anisotropic viscosity 
were obtained in \cite{LLT, ACW, CW2, DP1}.
For the 2D Boussinesq equations, the requirements on the initial data were significantly weakened in \cite{DP2, HKZ1,HKZ2}.  
Regarding the MHD-B\'enard system, 
some progress has been made in 2D case under various contexts, 
see, e.g., \cite{CD, DZ}. 
However, there has little work in the 3D case. 
Specifically, outstanding open problems such as global regularity of classic solutions 
for the fully dissipative system 
and whether the solutions blow up in finite time for the fully inviscid system 
remain unresolved. 

The main purpose of our paper is to obtain a Prodi-Serrin type regularity criterion 
for the non-diffusive 3D MHD-Boussinesq system. 
Unlike the case of the 3D Navier-Stokes equations, 
Prodi-Serrin type regularity criteria are not available for Euler equations in three-dimensional space.  
Thus, it is difficult to obtain global regularity for $u$, $b$, and $\theta$ simultaneously 
since there is no thermal diffusivity in the equation for $\theta$.
However, we are able to handle this by proving the higher order regularity for $u$ and $b$ first, 
before bounding $\Vert\nabla\theta\Vert_{L_{x}^2}$. We emphasize that this is the first work, to the best of our knowledge, that proves a Prodi-Serrin-type criterion in the case where the system is not fully dissipative.

The pioneering work of Serrin, Prodi, et al. (c.f. \cite{ESS, LM, So, St, Pr, S1, S2}) 
for the 3D Navier-Stokes equations proved that, for any $T>0$, if $u \in L_{t}^{r}([0,T];L_{x}^{s})$ 
with $2/r + 3/s < 1$ and $3 < s < \infty$, 
then the solution for the 3D Navier-Stokes equations remains regular on the interval $[0,T]$. 
Proof for the borderline case in various settings was obtained in \cite{ESS, LM, So, St}. 
Similar results concerning the 3D Navier-Stokes, Boussinesq and MHD equations were obtain in \cite{B, C, CT1, CT2, KZ1, KZ2, NP, P, PP1, PP2, PZ1, PZ2, S}. 
In particular, in \cite{Y1, Y2}, regularity criteria for MHD equations involving only two velocity components was proved but in a smaller Lebesgue space.
However, there is no literature on the regularity criteria for the solutions 
of systems (\ref{S2}) and (\ref{S1}). 
In this paper, we obtain a regularity condition using several velocity and magnetic field components 
in much larger space that is closer to the critical space of the equations. 
A central message of the present work is that 
with optimal and delicate application of our method, 
as well as potential new techniques such as in \cite{JZ, Q, S0, Y3, Ye1, Ye2}, 
one might further improve the criterion on the global regularity for system (\ref{S1}).

Moreover, we prove the local-in-time existence and uniqueness 
of the solutions to the system (\ref{S1}) with $H^3$ initial datum. 
We obtain the necessary a priori estimates 
and construct the solution via Galerkin methods 
for both the full and the non-diffusive systems. 
In particular, we show that the existence time of solutions to the full system 
does not depend on $\kappa$, 
which enables us to prove that the solutions to the full system approaches that of the non-diffusive system 
as $\kappa$ tends to $0$ on their time interval of existence. 
Regarding the fully inviscid system, we remark that the local well-posedness of 
either of the full system (\ref{S2}) or the non-diffusive system (\ref{S1})
is not automatically implied by that of the fully inviscid system (\ref{S3}), 
as observed in \cite{LT2} for multi-dimensional Burgers equation
\begin{align*}
     \frac{\partial u}{\partial t}
     +
     (u\cdot\nabla)u
     =
     \nu\Delta u,
\end{align*}
in two and higher dimensions. 
One might expect to that adding  more diffusion, namely in the form of a hyper-diffusion term $-\nu^2\Delta^2 u$, might make the equation even easier to handle.  However, the question well-posedness of the resulting equation, namely
\begin{align*}
    \frac{\partial u}{\partial t}
    +
    (u\cdot\nabla)u
    =
    -\nu^2\Delta^2 u
    +
    \nu\Delta u,
\end{align*}
remains open due to the lack of maximum principle, as observed in \cite{LT2}. 
Therefore, well-posedness is not automatic when additional diffusion is added, and it is worth exploring the regularity criteria of the solution to the non-diffusive and inviscid systems 
independent of the results for the full system. 
As we show in Section~\ref{sec3} and in the Appendix, we require a different approach to construct solutions, due to the lack of compactness in the non-dissipative system.
Note that the question of whether system (\ref{S3}) develops singularity in finite time still remains open. 

The paper is organized as follows.
In Section~\ref{sec2}, we provide the preliminaries for our subsequent work including the notation that we use, and state our main theorems.
In Section~\ref{sec3}, we obtain existence results for systems (\ref{S2}) and (\ref{S1}).
In Section~\ref{sec4}, we prove that solutions to the non-diffusive system (\ref{S1}) are unique. 
In Section~\ref{sec5}, we prove the regularity criteria for the solution to (\ref{S1}) using anisotropic estimates.
In the Appendix, we obtain the local in time well-posedness of the fully inviscid system (\ref{S3}) by a different argument.

\section{Preliminaries and summary of results}\label{sec2}
All through this paper we denote 
$\partial_{j} = \partial / \partial x_{j}$, 
$\partial_{jj} = \partial^{2}/\partial x_{j}^2$, 
$\partial_{t} = \partial /\partial t$, 
$\partial^{\alpha} = \partial^{|\alpha|} /\partial x_1^{\alpha_1}\cdots x_n^{\alpha_n}$, where $\alpha$ is a multi-index.  We also denote the horizontal gradient $\nabla_{h} = (\partial_1, \partial_2)$ 
and horizontal Laplacian $\Delta_{h} = \partial_{11} + \partial_{22}$. 
Also, we denote the usual Lebesgue and Sobolev spaces by $L_{x}^{p}$ and $H_{x}^{s}\equiv W_{x}^{s, 2}$, respectively, with the subscript $x$ (or $t$) indicating that the underlying variable is spatial (resp. temporal).  Let $\mathcal{F}$ be the set of all trigonometric polynomial over $\mathbb{T}^3$ 
and define the subset of divergence-free, zero-average trigonometric polynomials 
$$\mathcal{V} := \left\{ \phi\in\mathcal{F}: \nabla\cdot\phi = 0, \text{ and }\int_{\mathbb{T}^3}\phi\,dx= 0\right\}.$$
We use the standard convention of denoting by $H$ and $V$ the closures of $\mathcal{V}$ in $L_{x}^2$ and $H_{x}^1$, respectively, 
with inner products 
$$(u, v) = \sum_{i=1}^3\int_{\mathbb{T}^3}u_{i}v_{i}\,dx \text{ \,\,and\,\,  } (\nabla u, \nabla v) = \sum_{i, j=1}^3\int_{\mathbb{T}^3}\partial_{j}u_{i}\partial_{j}v_{i}\,dx,$$
respectively, associated with the norms $|u|=(u, u)^{1/2}$ and $\Vert u \Vert=(\nabla u, \nabla u)^{1/2}$.
The latter is a norm due to the Poincar\'e inequality $$\Vert\phi\Vert_{L_{x}^2} \leq C\Vert\nabla\phi\Vert_{L_{x}^2}$$ 
holding for all  $\phi\in V$.
We also have the following compact embeddings (see, e.g., \cite{CF,T1})
$$V \hookrightarrow H \hookrightarrow V',$$
where $V'$ denotes the dual space of $V$.

The following interpolation result is frequently used in this paper (see, e.g., \cite{N} for a detailed proof).
Assume $1 \leq q, r \leq \infty$, and $0<\gamma<1$.  
For $v\in L_{x}^q(\mathbb{T}^{n})$, such that  $\partial^\alpha v\in L_{x}^{r} (\mathbb{T}^{n})$, for $|\alpha|=m$, then 
\begin{align}\label{PT1}
\Vert\partial_{s}v\Vert_{L^{p}} \leq C\Vert\partial^{\alpha}v\Vert_{L^{r}}^{\gamma}\Vert v\Vert_{L^{q}}^{1-\gamma},
\quad\text{where}\quad
\frac{1}{p} - \frac{s}{n} = \left(\frac{1}{r} - \frac{m}{n}\right) \gamma+ \frac{1}{q}(1-\gamma).
\end{align}

The following materials are standard in the study of fluid dynamics, 
in particular for the Navier-Stokes equations,  
and we refer to reader to \cite{CF, T1} for more details. 
We define the Stokes operator $A\triangleq -P_{\sigma}\Delta$ 
with domain $\mathcal{D}(A)\triangleq H_{x}^2\cap V$, 
where $P_{\sigma}$ is the Leray-Helmholtz projection. 
Note that under periodic boundary conditions, 
we have $A = -\Delta P_{\sigma}$.
Moreover, the Stokes operator can be extended 
as a linear operator from $V$ to $V'$ as 
$$\left<Au, v\right> = (\nabla u, \nabla v) \text{  for all  } v\in V.$$
It is well-known that $A^{-1} : H \hookrightarrow \mathcal{D}(A)$ 
is a positive-definite, self-adjoint, and compact operator from $H$ into itself,
thus, $A^{-1}$ possesses an orthonormal basis of positive eigenfunctions $\{ w_{k}\}_{k=1}^{\infty}$ 
in $H$, corresponding to a sequence of non-increasing sequence of eigenvalues. 
Therefore, $A$ has non-decreasing eigenvalues $\lambda_{k}$, 
i.e., $0 \leq \lambda_1 \leq \lambda_2, \ldots$ 
since $\{ w_{k}\}_{k=1}^{\infty}$ are also eigenfunctions of $A$. 
Furthermore, for any integer $M > 0$, we define $H_{M}\triangleq \text{span}\{w_1, w_2, \ldots, w_{M}\}$ 
and $P_{M} : H \to H_{M}$ be the $L_{x}^2$ orthogonal projection onto $H_{M}$. 
Next, for any $u, v, w \in \mathcal{V}$, 
we introduce the convenient notation for the bilinear term
\begin{align*}
     B(u, v) := P_{\sigma}((u\cdot\nabla)v),
\end{align*}
which can be extended to a continuous map 
$B : V \times V \to V'$ such that
\begin{align*}
     \left<B(u, v), w\right> = \int_{\mathbb{T}^3}(u\cdot\nabla v)\cdot w\,dx.
\end{align*}
for smooth functions $u,v,w\in V$.  
Notice that $\theta$ is a scalar function 
so we cannot actually apply $P_{\sigma}$ on it; 
hence, the notation $P_{M}\theta$ 
should be understood as projection onto the space spanned by the first $M$ eigenfunctions of $-\Delta$ only. 
Therefore, in order to avoid abuse of notation, we denote $\mathcal{B}(u, \theta):= u\cdot\nabla\theta$ for smooth functions, and extended it to a continuous map $B : V \times H^1 \to H^{-1}$ similarly to $B(\cdot,\cdot)$. 
We will use the following important properties of the map $B$. Detailed proof can be found in, e.g., \cite{CF, FMRT}.
\begin{Lemma}
\label{PL1}
For the operator $B$, we have 
\begin{subequations}
\begin{align}
    \label{embd1}
    \ip{B(u,v)}{w}_{V'} &= -\ip{B(u,w)}{v}_{V'}, 
    &\quad\forall\;u\in V, v\in V, w\in V,\\
    \label{embd2}
    \ip{B(u,v)}{v}_{V'} &= 0,
    &\quad\forall\;u\in V, v\in V, w\in V,\\
    \label{B:326}
    |\ip{B(u,v)}{w}_{V'}|
    &\leq C\Vert u\Vert_{L_{x}^2}^{1/2} \Vert\nabla u\Vert_{L_{x}^2}^{1/2} \Vert\nabla v\Vert_{L_{x}^2} \Vert\nabla w\Vert_{L_{x}^2},
    &\quad\forall\;u\in V, v\in V, w\in V,\\
    \label{B:623}
    |\ip{B(u,v)}{w}_{V'}|
    &\leq C\Vert\nabla u\Vert_{L_{x}^2} \Vert\nabla v\Vert_{L_{x}^2} \Vert w\Vert_{L_{x}^2}^{1/2} \Vert\nabla w\Vert_{L_{x}^2}^{1/2},
    &\quad\forall\;u\in V, v\in V, w\in V,\\
    \label{B:236}
    |\ip{B(u,v)}{w}_{V'}|
    &\leq C\Vert u\Vert_{L_{x}^2} \Vert\nabla v\Vert_{L_{x}^2}^{1/2} \Vert Av\Vert_{L_{x}^2}^{1/2} \Vert\nabla w\Vert_{L_{x}^2},
    &\quad\forall\;u\in H, v\in \mathcal{D}(A), w\in V,\\
    \label{B:632}
    |\ip{B(u,v)}{w}_{V'}|
    &\leq C\Vert\nabla u\Vert_{L_{x}^2} \Vert\nabla v\Vert_{L_{x}^2}^{1/2} \Vert Av\Vert_{L_{x}^2}^{1/2} \Vert w\Vert_{L_{x}^2},
    &\quad\forall\;u\in V, v\in \mathcal{D}(A), w\in H,\\
    \label{B:i22}
    |\ip{B(u,v)}{w}_{V'}|
    &\leq C\Vert\nabla u\Vert_{L_{x}^2}^{1/2} \Vert Au\Vert_{L_{x}^2}^{1/2} \Vert\nabla v\Vert_{L_{x}^2} \Vert w\Vert_{L_{x}^2},
    &\quad\forall\;u\in \mathcal{D}(A), v\in V, w\in H,\\
    \label{B:263}
    |\ip{B(u,v)}{w}_{V'}|
    &\leq C\Vert u\Vert_{L_{x}^2} \Vert Av\Vert_{L_{x}^2} \Vert w\Vert_{L_{x}^2}^{1/2} \Vert\nabla w\Vert_{L_{x}^2}^{1/2},
    &\quad\forall\;u\in H, v\in \mathcal{D}(A), w\in V,\\
    \label{B:362s}
    |\ip{B(u,v)}{w}_{\mathcal{D}(A)'}|
    &\leq C\Vert u\Vert_{L_{x}^2}^{1/2} \Vert\nabla u\Vert_{L_{x}^2}^{1/2} \Vert v\Vert_{L_{x}^2} \Vert Aw\Vert_{L_{x}^2},
    &\quad \forall\;u\in V, v\in H, w\in\mathcal{D}(A).
\end{align}
\end{subequations}

Moreover, essentially identical results hold for $\mathcal{B}(u, \theta)$, \emph{mutatis mutandis}.

\end{Lemma}

The following lemma is a special case of the Troisi inequality from \cite{T} and is useful for our estimates throughout the paper. 
\begin{Lemma}
\label{PL2}
There exists a constant $C>0$ such that for $v \in C_{0}^{\infty}(\mathbb{R}^3)$, we have
$$\Vert v \Vert_{L^{6}} \leq C\prod_{i=1}^{3}\Vert\partial_{i}v\Vert_{L^2}^{\frac{1}{3}}.$$
\end{Lemma}

Regarding the pressure term, we recall the fact that, for any distribution $f$, 
the equality $f=\nabla p$ holds for some distribution $p$ if and only if $\left<f, w\right> = 0$ for all $w \in \mathcal {V}$. 
See \cite{W1} for details.

Next, we list three fundamental lemmas needed in order to prove Theorem~\ref{T0}.  Their proofs can be found in \cite{KZ2} and \cite{Y2}, respectively.
\begin{Lemma}
\label{L01}
Assume $u = (u_1, u_2, u_3) \in H^2(\mathbb{T}^3)\cap V$. Then
\begin{align*}
     \sum_{j, k=1}^2 \int_{\mathbb{T}^3} u_{j}\partial_{j}u_{k} \Delta_{h}u_{k}\,dx
     =
     \frac{1}{2}\sum_{j, k=1}^2 \int_{\mathbb{T}^3} \partial_{j}u_{k}\partial_{j}u_{k}\partial_{3}u_{3}\,dx
     -
     \int_{\mathbb{T}^3} \partial_{1}u_{1}\partial_{2}u_{2}\partial_{3}u_{3}\,dx
     +
     \int_{\mathbb{T}^3} \partial_{1}u_{2}\partial_{2}u_{1}\partial_{3}u_{3}\,dx.
\end{align*}
\end{Lemma}

\begin{Lemma}
\label{L02}
For $u$ and $b$ from the solution of (\ref{S1}) and $i = 1, 2, 3$, we have
\begin{align*}
     &\int_{\mathbb{T}^3} u_{j}\partial_{j}u_{k} \partial_{ii}u_{k}\,dx
     -
     \int_{\mathbb{T}^3} b_{j}\partial_{j}b_{k} \partial_{ii}u_{k}\,dx
     +
     \int_{\mathbb{T}^3} u_{j}\partial_{j}b_{k} \partial_{ii}b_{k}\,dx
     -
     \int_{\mathbb{T}^3} b_{j}\partial_{j}u_{k} \partial_{ii}b_{k}\,dx
     \\&
     =
     \sum_{j, k = 1}^{3}\int_{\mathbb{T}^3} -\partial_{i}u_{j}\partial_{j}u_{k}\partial_{i}u_{k}\,dx
     +
     \int_{\mathbb{T}^3} \partial_{i}b_{j}\partial_{j}b_{k}\partial_{i}u_{k}\,dx
     -
     \int_{\mathbb{T}^3} \partial_{i}u_{j}\partial_{j}b_{k}\partial_{i}b_{k}\,dx
     +
     \int_{\mathbb{T}^3} \partial_{i}b_{j}\partial_{j}u_{k}\partial_{i}b_{k}\,dx.
\end{align*}
\end{Lemma}

The following Aubin-Lions Compactness Lemma is needed in order to construct solutions for (\ref{S2}).
\begin{Lemma}
\label{ALT}
Let $T>0$, $p \in (1, \infty)$ and let $\{f_{n}(t, \cdot)\}_{n=1}^{\infty}$ 
be a bounded sequence of function in $L_{t}^{p}([0, T]; Y)$ 
where $Y$ is a Banach space. If $\{f_{n}\}_{n=1}^{\infty}$ 
is also bounded in $L_{t}^{p}([0, T]; X)$, where $X$ is compactly imbedded in $Y$ 
and $\{\partial f_{n}/\partial t\}_{n=1}^{\infty}$ is bounded in $L_{t}^{p}([0, T]; Z)$ uniformly 
where $Y$ is continuously imbedded in $Z$, 
then $\{f_{n}\}_{n=1}^{\infty}$ is relatively compact in $L_{t}^{p}([0, T]; Y)$.
\end{Lemma}

The following theorem is our main result.
 It provides a Prodi-Serrin type regularity criterion for system (\ref{S1}).
\begin{Theorem}
\label{T0}
For $m\geq3$, $u_0$, $b_0 \in H_{x}^{m}\cap V$, 
and $\theta_0 \in H_{x}^3$, 
if we further assume that $u_2, u_3, b_2, b_3 \in L_{t}^{r}([0, T); L_{x}^{s}(\mathbb{T}^{3}))$ and 
\begin{align*}
     &\frac{2}{r} + \frac{3}{s} 
     =
     \frac{3}{4}
     +
     \frac{1}{2s}, 
     \qquad s > 10/3,
\end{align*}
for a given $T > T^{*}$ where $T^{*}$ is in Theorem~\ref{T3}.
Then $(u, b, \theta)$ remains smooth beyond $T^{*}$. 
Namely, $\Vert u\Vert_{H_{x}^1}$, $\Vert b\Vert_{H_{x}^1}$, and $\Vert\theta\Vert_{H_{x}^1}$ 
remains bounded up to $T > T^{*}$, 
and consequently, $u, b, \theta \in C^{\infty}(\Omega \times (0, T))$.
\end{Theorem}

The next three theorems provide local well-posedness for systems (\ref{S3}), (\ref{S2}), and (\ref{S1}). 
First, for the fully inviscid system (\ref{S3}), we have
\begin{Theorem}
\label{T1}
For the initial data $(u_0, b_0, \theta_0) \in H_{x}^{3}\cap V$, 
there exists a unique solution 
$$(u, b, \theta) \in L_{t}^{\infty}((0, \tilde{T}); H_{x}^{3}\cap V)$$ 
to the fully inviscid MHD-Boussinesq system (\ref{S3}) 
for some $\tilde{T} > 0$, depending on $g$ and the initial data.
\end{Theorem}

Regarding system (\ref{S2}), we have
\begin{Theorem}
\label{T2}
For $m\geq3$ and $u_0, b_0 \in H_{x}^{m} \cap V$, 
and $\theta_0 \in H_{x}^{m}$, 
there exists a solution $(u, b, \theta)$ with $u, b \in C_{w}([0, T); H) \cap L_{t}^{2}((0, T); V)$ 
and $\theta \in C_{w}([0, T); L_{x}^2) \cap L_{t}^{2}((0, T); H_{x}^1)$ for any $T > 0$ 
for (\ref{S2}).
Also, the solution is unique if $u, b \in L_{t}^{\infty}([0, T'); H_{x}^{m}\cap V) \cap L_{t}^{2}((0, T'); H_{x}^{m+1} \cap V)$ 
and $\theta \in L_{t}^{\infty}([0, T'); H_{x}^{m}) \cap L_{t}^{2}((0, T'); H_{x}^{m+1})$ 
with some $T'$ depending only on $\nu$, $\eta$, and the initial datum. 
\end{Theorem}

For the non-diffusive MHD-Boussinesq system (\ref{S1}), which we mainly focus on, we have
\begin{Theorem}
\label{T3}
For $m\geq3$ and $u_0$, $b_0 \in H_{x}^{m} \cap V$, $\theta_0 \in H_{x}^{m}$, 
there exists a unique solution $(u, b, \theta)$ to the non-diffusive MHD-Boussinesq system (\ref{S1}),
where $u, b \in L_{t}^{\infty}([0, T^{*}); H_{x}^{m}\cap V) \cap L_{t}^{2}((0, T^{*}); H_{x}^{m+1}\cap V)$ divergence free, and $\theta \in L_{t}^{\infty}([0, T^{*}); H_{x}^{m})$, 
where $T^{*}$ depends on $\nu$, $\eta$, and the initial datum. 
\end{Theorem}

\section{Proof of Theorem~\ref{T2} and Theorem~\ref{T3}}\label{sec3}
For Theorem~\ref{T2}, 
we use Galerkin approximation to obtain the solution for the full MHD-Boussinesq system (\ref{S2}), 
while for the existence part of Theorem~\ref{T3}, 
the proof is similar with only minor modification so we omit the details.

{\smallskip\noindent {\em Proof of Theorem~\ref{T2}.}}
Consider the following finite-dimensional ODE system, which we think of as an approximation to system \eqref{S2} after apply the Leray projection $P_\sigma$.  
\begin{equation}
  \left\{
     \begin{aligned}
     &
     \frac{d u_{M}}{dt}
     -
     \nu Au_{M}
     +
     P_{M}B(u_{M}, u_{M})
     =
     P_{M}B(b_{M}, b_{M})
     +
     gP_{\sigma}(\theta_{M}e_{3}),
     \\&
     \frac{d b_{M}}{dt}
     -
     \eta Ab_{M}
     +
     P_{M}B(u_{M}, b_{M})
     =
     P_{M}B(b_{M}, u_{M}),
     \\&
     \frac{d \theta_{M}}{dt}
     -
     \kappa\Delta\theta_{M}
     +
     P_{M}\mathcal{B}(u_{M}, \theta_{M})
     =
     0,
    \end{aligned}
  \right.
  \label{S4}
\end{equation}
with initial datum $P_{M}u(\cdot, 0)=u_{M}(0)$, $P_{M}b(\cdot, 0)=b_{M}(0)$, and $P_{M}\theta(\cdot, 0)=\theta_{M}(0)$. 
Notice that all terms but the time-derivatives of the above ODE systems are at most quadratic, and therefore they are locally Lipschitz continuous.  Thus, by the Picard-Lindelhoff Theorem, we know that there exists a solution up to some time $T_{M} > 0$. 
Next we take justified inner-products with the above three equations by $u_{M}$, $b_{M}$, and $\theta_{M}$, respectively, integrate by parts, and add the results to obtain
\begin{align*}
     &\quad
     \frac{1}{2}\frac{d}{dt}
     \left( \Vert u_{M}\Vert_{L_{x}^2}^2 + \Vert b_{M}\Vert_{L_{x}^2}^2 + \Vert \theta_{M}\Vert_{L_{x}^2}^2 \right)
     +
     \nu\Vert\nabla u_{M}\Vert_{L_{x}^2}^2
     +
     \eta\Vert\nabla b_{M}\Vert_{L_{x}^2}^2
     +
     \kappa\Vert\nabla \theta_{M}\Vert_{L_{x}^2}^2
     \\&
     =
     \int_{\mathbb{T}^3}(b_{M}\cdot\nabla)b_{M}u_{M}\,dx
     +
     \int_{\mathbb{T}^3}g\theta_{M}u_{M}e_{3}\,dx
     +
     \int_{\mathbb{T}^3}(b_{M}\cdot\nabla)u_{M}b_{M}\,dx
     \\&
     =
     g\int_{\mathbb{T}^3}\theta_{M}u_{M}e_{3}\,dx,    
\end{align*}
where we used the divergence free condition, Lemma~\ref{PL1}, and the orthogonality of $P_{\sigma}$ and $P_{M}$.  By the Cauchy-Schwarz and Young's inequalities, 
we obtain
\begin{align}
     &\quad
     \frac{d}{dt}
     \left( \Vert u_{M}\Vert_{L_{x}^2}^2 + \Vert b_{M}\Vert_{L_{x}^2}^2 + \Vert \theta_{M}\Vert_{L_{x}^2}^2 \right)
     +
     2\nu\Vert\nabla u_{M}\Vert_{L_{x}^2}^2
     +
     2\eta\Vert\nabla b_{M}\Vert_{L_{x}^2}^2
     +
     2\kappa\Vert\nabla \theta_{M}\Vert_{L_{x}^2}^2 \nonumber
     \\&
     \leq
     C_{g} \left( \Vert u_{M}\Vert_{L_{x}^2}^2 + \Vert \theta_{M}\Vert_{L_{x}^2}^2\right).
     \label{Ineq1}
\end{align}
Thus, by the differential form of Gr\"onwall's inequality,  $u_{M}$ and $b_{M}$ are uniformly bounded in $L_{t}^{\infty}([0, T_M); H)$, 
while $\theta_{M}$ is uniformly bounded in $L_{t}^{\infty}([0, T_M); L_{x}^2$, independently of $T_M$. 
Namely, 
\begin{align*}
     \Vert u_{M}(t)\Vert_{L_{x}^2}^2 + \Vert b_{M}(t)\Vert_{L_{x}^2}^2 + \Vert \theta_{M}(t)\Vert_{L_{x}^2}^2
     \leq
     C_{g, T} \Vert u_{M}(0)\Vert_{L_{x}^2}^2 + \Vert b_{M}(0)\Vert_{L_{x}^2}^2 + \Vert \theta_{M}(0)\Vert_{L_{x}^2}^2,
\end{align*}
for any $0<t<T_M$.  Thus, for each $M$, the solutions can be extended uniquely beyond $T_M$ to an interval $[0,T]$, where $T>0$ is arbitrary.  In particular, the interval of existence and uniqueness is independent of $M$.  Using the embedding $L_{t}^{\infty}\hookrightarrow L_{t}^{2}$, and extracting a subsequence if necessary (which we relabel as $(u_M, b_M, \theta_M)$), we may invoke the Banach-Alaoglu Theorem to obtain $u, b\in L_{t}^{2}([0, T]; H)$, and $\theta \in L_{t}^2([0, T]; L_{x}^2)$, such that
\begin{align*}
     &
     u_{M} \rightharpoonup u \text{ \,\,\,and\,\,\, } b_{M} \rightharpoonup b \text{ \,\,\,weakly in \,\,\,} L_{t}^2([0, T];H),
     \\&
     \theta_{M} \rightharpoonup \theta \text{ \,\,\,weakly in\,\,\, } L_{t}^2([0, T];L_{x}^2).
\end{align*}
$(u,b,\theta)$ is our candidate solution.  Next, integrating (\ref{Ineq1}) over time from $0$ to $t < T$, and using Gr\"onwall's inequality,  we have that $u_{M}$ and $b_{M}$ are uniformly bounded in $L_{t}^{2}([0, t); V)$, while $\theta_{M}$ is uniformly bounded in $L_{t}^{2}([0, T); H_{x}^{1})$ for any $T>0$.  Next, we obtain bounds on 
$d u_{M}/dt$, $d b_{M}/dt$, and $d \theta_{M}/dt$ 
in certain functional space uniformly with respect to $M$.
Note that
\begin{equation}
  \left\{
     \begin{aligned}
     &
     \frac{d u_{M}}{dt}
     =
     -\nu Au_{M}
     -   
     P_{M}B(u_{M}, u_{M})
     +
     P_{M}B(b_{M}, b_{M})
     +
     gP_{M}(\theta_{M}e_{3}),
     \\&
     \frac{d b_{M}}{dt}
     =
     -\eta Ab_{M}
     -
     P_{M}B(u_{M}, b_{M})
     +
     P_{M}B(b_{M}, u_{M}),
     \\&
     \frac{d \theta_{M}}{dt}
     =
     -\kappa\Delta\theta_{M}
     -
     \mathcal{B}((u_{M}, \theta_{M}).
    \end{aligned}
 \right.
\end{equation}
Note in the first equation that $Au_{M}$ is bounded in $L_{t}^{2}([0, T); V')$ 
due to the fact that $u_{M}$ is bounded in $L_{t}^{2}([0, T); V)$. 
Also, we have $gP_{M}(\theta_{M}e_{3})$ is bounded in $L_{t}^{2}([0, T); H)$. 
On the other hand, by Lemma~\ref{PL1}, 
we have $$\Vert P_{M}B(u_{M}, u_{M})\Vert_{V'} \leq C\Vert u_{M}\Vert_{L_{x}^2}^{1/2}\Vert\nabla u_{M}\Vert_{L_{x}^2}^{3/2},$$
as well as $$\Vert P_{M}B(b_{M}, b_{M})\Vert_{V'} \leq C\Vert b_{M}\Vert_{L_{x}^2}^{1/2}\Vert\nabla b_{M}\Vert_{L_{x}^2}^{3/2}.$$ 
Since the $L^2$-norm of $u_{M}$ is uniformly bounded 
and the $L^2$-norm of $\nabla u_{M}$ are uniformly integrable, 
we see that $d u_{M}/dt$ are bounded in $L_{t}^{4/3}([0, T); V')$. 
Similarly, from the second and third equations, 
we have that $d b_{M}/dt$ and $d \theta_{M}/dt$ 
are also bounded in $L_{t}^{4/3}([0, T); V')$ 
and $L_{t}^{4/3}([0, T); H_{x}^{-1})$, respectively. 
Therefore, by Lemma~\ref{ALT} and the uniform bounds obtained above, 
there exists a subsequence (which we again relabel as $(u_{M}, b_{M}, \theta_{M})$ if necessary)  such that 
\begin{align*}
     &
     u_{M} \to u \text{ \,\,\,and\,\,\, } b_{M} \to b \text{ \,\,\,strongly in \,\,\,} L_{t}^2([0, T];H),
     \\&
     \theta_{M} \to \theta \text{ \,\,\,strongly in\,\,\, } L_{t}^2([0, T];L_{x}^2), 
     \\&
     u_{M} \rightharpoonup u \text{ \,\,\,and\,\,\, } b_{M} \to b \text{ \,\,\,weakly in\,\,\, } L_{t}^2([0, T];V),
     \\&
     \theta_{M} \rightharpoonup \theta \text{ \,\,\,weakly in\,\,\, } L_{t}^2([0, T];H_{x}^1), 
     \\&
     u_{M} \rightharpoonup u \text{ \,\,\,and\,\,\, } b_{M} \to b \text{ \,\,\,weak-$\ast$ in\,\,\, } L_{t}^{\infty}([0, T];H),
     \\&
     \theta_{M} \rightharpoonup \theta \text{ \,\,\,weak-$\ast$ in\,\,\, } L_{t}^{\infty}([0, T];L_{x}^2),
\end{align*}
for any $T>0$.
Thus, by taking inner products of (\ref{S4}) with test function $\psi(t, x) \in C_{t}^1([0, T]; C_{x}^{\infty})$ 
with $\psi(T) = 0$, and using the standard arguments of strong/weak convergence for Navier-Stokes equations (see, e.g., \cite{CF,T1}), 
we have that each of the linear and nonlinear terms in (\ref{S4}) converges to the appropriate limit in an appropriate weak sense. 
Namely, we obtain that (\ref{S2}) holds in the weak sense, 
where the pressure term $p$ is recovered by the approach mentioned in Section~\ref{sec2} and we omit the details here. 
Finally, we take action of (\ref{S2}) with an arbitrary $v \in \mathcal{V}$. 
Then, by integrating in time over $[t_0, t_1] \subset [0, T]$ and sending $t_1 \to t_0$ 
one can prove by standard arguments (c.f. \cite{CF,T1}) that $u, b$ and $\theta$ are in fact weakly continuous in time. 
Therefore, the initial condition is satisfied in the weak sense.

Next we show that the solution is in fact regular at least for short time, provided $(u_0, b_0, \theta_0) \in H^{m}\cap V$. 
We start by multiplying (\ref{S2}) by $ Au$, $ Ab$, and $\Delta \theta$, respectively, 
integrate over $\mathbb{T}^3$, and add, to obtain
\begin{align*}
     &\quad
     \frac{1}{2}\frac{d}{dt}
     \left( \Vert \nabla u\Vert_{L_{x}^2}^2 + \Vert \nabla b\Vert_{L_{x}^2}^2 + \Vert \nabla \theta\Vert_{L_{x}^2}^2 \right)
     +
     \nu\Vert\Delta u\Vert_{L_{x}^2}^2
     +
     \eta\Vert\Delta b\Vert_{L_{x}^2}^2
     +
     \kappa\Vert\Delta \theta\Vert_{L_{x}^2}^2
     \\&
     =
     -\int_{\mathbb{T}^3}(u\cdot\nabla)u\Delta u\,dx
     +
     \int_{\mathbb{T}^3}(b\cdot\nabla)b\Delta u\,dx
     +
     g\int_{\mathbb{T}^3}\theta\Delta ue_{3}\,dx, 
     \\&
     \quad-\int_{\mathbb{T}^3}(u\cdot\nabla)b\Delta b\,dx
     +
     \int_{\mathbb{T}^3}(b\cdot\nabla)u\Delta b\,dx
     -
     \int_{\mathbb{T}^3}(u\cdot\nabla)\theta\Delta \theta\,dx
     \\&
     \leq
     C\Vert\nabla u\Vert_{L_{x}^2}^{3/2}\Vert\Delta u\Vert_{L_{x}^2}^{3/2}
     +
     C\Vert\nabla b\Vert_{L_{x}^2}^{3/2}\Vert\Delta b\Vert_{L_{x}^2}^{1/2}\Vert\Delta u\Vert_{L_{x}^2}
     +
     g\Vert\nabla u\Vert_{L_{x}^2}\Vert\nabla\theta\Vert_{L_{x}^2}
     \\&\quad
     +
     C\Vert\nabla u\Vert_{L_{x}^2}\Vert\nabla b\Vert_{L_{x}^2}^{1/2}\Vert\Delta b\Vert_{L_{x}^2}^{3/2}
     +
     C\Vert\nabla b\Vert_{L_{x}^2}\Vert\nabla u\Vert_{L_{x}^2}^{1/2}\Vert\Delta u\Vert_{L_{x}^2}^{1/2}\Vert\Delta b\Vert_{L_{x}^2}
     \\&\quad
     +
     C\Vert\theta\Vert_{L_{x}^{\infty}}\Vert\nabla u\Vert_{L_{x}^2}\Vert\Delta\theta\Vert_{L_{x}^2}
     \\&
     \leq
     \frac{\nu}{2}\Vert\Delta u\Vert_{L_{x}^2}^{2}
     +
     \frac{\eta}{2}\Vert\Delta b\Vert_{L_{x}^2}^{2}
     +
     \frac{\kappa}{2}\Vert\Delta\theta\Vert_{L_{x}^2}^{2}
     \\&\quad
     +
     \frac{C}{\nu^3}\Vert\nabla u\Vert_{L_{x}^2}^{6}
     +
     \frac{C}{\nu\eta}\Vert\nabla b\Vert_{L_{x}^2}^{6}
     +
     C\Vert\nabla\theta\Vert_{L_{x}^2}^{2}
     +
     C\Vert\nabla u\Vert_{L_{x}^2}^{2}
     \\&\quad
     +
     \frac{C}{\eta^3}\Vert\nabla u\Vert_{L_{x}^2}^{4}\Vert\nabla b\Vert_{L_{x}^2}^{2}
     +
     \frac{C}{\nu\eta}\Vert\nabla b\Vert_{L_{x}^2}^{4}\Vert\nabla u\Vert_{L_{x}^2}^{2}
     +
     \frac{C}{\kappa}\Vert\nabla u\Vert_{L_{x}^2}^{2},
\end{align*}
where we applied the H\"older's inequality, Sobolev embedding, and Young's inequality.
By denoting 
$$K(t) = \Vert \nabla u(t)\Vert_{L_{x}^2}^2 + \Vert \nabla b(t)\Vert_{L_{x}^2}^2 + \Vert \nabla \theta(t)\Vert_{L_{x}^2}^2,$$ 
we have $$\frac{d K}{d t} \leq CK + CK^3,$$
which implies that there exists a $T'>0$ such that 
\begin{equation}\label{H1_bound_K1}
     K(t)
     \leq
     \frac{Ce^{CT'/2}K(0)}{\sqrt{1-K^2(0)(e^{CT'} - 1)}}=:K_1(T'),\quad\text{for all } t\in[0,T'].
\end{equation}
After integrating from $t=0$ to $t=T'$ and the constant $C$ depends on the initial datum, $g$, $\nu$, $\eta$, and $\kappa$. 
This shows that $(u, b, \theta) \in L_{t}^{\infty}((0, T'); H^1\cap V)$ as $M \to \infty$, provided $T' < 1/K^2(0)e^{2C}$.

In order to pass to the limit $\kappa \to 0^{+}$, 
we must show that the above existence time $T'$ is independent of $\kappa$. 
We follow the vanishing viscosity technique for the Navier-Stokes equations, (c.f. \cite{CF})
i.e., let $\tau = \kappa t$, 
and denote
\begin{align*}
     \tilde{Q}(\tau) 
     =
     \frac{1}{\kappa}
     \left(
     \Vert\nabla u(\frac{\tau}{\kappa})\Vert_{L_{x}^2}
     +
     \Vert\nabla b(\frac{\tau}{\kappa})\Vert_{L_{x}^2}
     +
     \Vert\nabla\theta(\frac{\tau}{\kappa})\Vert_{L_{x}^2}
     \right).
\end{align*}
The above $H^1$ estimates thus imply that
\begin{align*}
     \frac{d\tilde{Q}}{d\tau}
     \leq
     \tilde{C} + \tilde{C}\tilde{Q}^2,
\end{align*}
where $\tilde{C}$ depends only on $g$, $\nu$, $\eta$, and is independent of $\kappa$.
Thus, integrating from $\tau=0$ to $\tau=\tilde{\tau}$,
we obtain 
\begin{align*}
     \tilde{Q}(\tilde{\tau})
     \leq
     \frac{\tilde{Q}(0)}{1 - \tilde{C}\tilde{\tau}\tilde{Q}(0)}.
\end{align*}
Thus, if $$\tilde{C}\tilde{\tau}\tilde{Q}(0) \leq \delta < 1, $$ i.e., 
\begin{align*}
     \tilde{C} (\kappa\tilde{t}) \frac{1}{\kappa}\left(
     \Vert\nabla u(0)\Vert_{L_{x}^2}
     +
     \Vert\nabla b(0)\Vert_{L_{x}^2}
     +
     \Vert\nabla\theta(0)\Vert_{L_{x}^2}\right)
     \leq
     \delta < 1,
\end{align*}
it follows that $\tilde{Q}(\tilde{\tau}) \leq C_{\delta}\tilde{Q}(0)$.
Hence, we have proved that, if
\begin{equation}
     T'
     <
     \frac{\tilde{C}}
     {\left(
     \Vert\nabla u(0)\Vert_{L_{x}^2}
     +
     \Vert\nabla b(0)\Vert_{L_{x}^2}
     +
     \Vert\nabla\theta(0)\Vert_{L_{x}^2}\right)},
  \label{Ineq_time}
\end{equation}
then the above $H^1$ estimates remain valid for any $\kappa > 0$.

On the other hand, we showed earlier that 
$$\nu\int_{0}^{T'}\Vert\Delta u\Vert_{L_{x}^2}^2\,dt + \eta\int_{0}^{T'}\Vert\Delta b\Vert_{L_{x}^2}^2\,dt + \kappa\int_{0}^{T'}\Vert\Delta\theta\Vert_{L_{x}^2}^2\,dt$$
remains bounded as $M \to \infty$.
Thus, we have $(u, b, \theta) \in L_{t}^2((0, T'); H^2\cap V)$. 
In order to obtain the higher-order regularity in $H^2$ and $H^3$, 
we follow standard arguments (see, e.g., \cite{MP}) and apply the following argument successively. 
First, for a multi-index $\alpha$ of order $|\alpha|=2$,  
we apply the partial differential operator $\partial^{\alpha}$, 
to (\ref{S2}), 
and test the equations for $u$, $b$, and $\theta$ 
by $\partial^{\alpha} u$, $\partial^{\alpha} b$, and $\partial^{\alpha} \theta$, respectively, and obtain
\begin{align*}
  \left\{
     \begin{aligned}
     &
     \frac{1}{2}\frac{d}{dt}\Vert\partial^{\alpha} u\Vert_{L_{x}^2}^2
     +
     \nu\Vert\nabla\partial^{\alpha} u\Vert_{L_{x}^2}^2
     &&=
     \int_{\mathbb{T}^3}\partial^{\alpha}((b\cdot\nabla) b)\partial^{\alpha} u\,dx
     -
     \int_{\mathbb{T}^3}\partial^{\alpha}((u\cdot\nabla) u)\partial^{\alpha} u\,dx
     +
     g\int_{\mathbb{T}^3}\partial^{\alpha}\theta \partial^{\alpha} u\,dx
     \\&
     &&=
     I_1 + I_2 + I_3,
     \\&
     \frac{1}{2}\frac{d}{dt}\Vert\partial^{\alpha} b\Vert_{L_{x}^2}^2
     +
     \eta\Vert\nabla\partial^{\alpha} b\Vert_{L_{x}^2}^2
     &&=
     \int_{\mathbb{T}^3}\partial^{\alpha}((b\cdot\nabla) u)\partial^{\alpha} b\,dx
     -
     \int_{\mathbb{T}^3}\partial^{\alpha}((u\cdot\nabla) b)\partial^{\alpha} b\,dx
     =
     I_4 + I_5,
     \\&
     \frac{1}{2}\frac{d}{dt}\Vert\partial^{\alpha} \theta\Vert_{L_{x}^2}^2
     +
     \kappa\Vert\nabla\partial^{\alpha} \theta\Vert_{L_{x}^2}^2
     &&=
     -\int_{\mathbb{T}^3}\partial^{\alpha}((u\cdot\nabla) \theta)\partial^{\alpha} \theta\,dx
     =
     I_6.
   \end{aligned}
 \right.
\end{align*}
In order to estimate $I_1$, we use Lemma~\ref{PL1} and get
\begin{align*}
     I_1
     &
     =
     \sum_{\zeta\leq\alpha}\binom{\alpha}{\zeta}\int_{\mathbb{T}^3}((\partial^{\zeta}b\cdot\nabla) \partial^{\alpha-\zeta}b)\partial^{\alpha} u\,dx
     \\&
     \leq
     C\Vert\nabla b\Vert_{L_{x}^2}\Vert\partial^{\alpha} u\Vert_{L_{x}^2}^{1/2}\Vert\nabla\partial^{\alpha} u\Vert_{L_{x}^2}^{1/2}\Vert\nabla\partial^{\alpha} b\Vert_{L_{x}^2}
     +
     C\Vert\nabla b\Vert_{L_{x}^2}\Vert\partial^{\alpha} u\Vert_{L_{x}^2}^{1/2}\Vert\nabla\partial^{\alpha} u\Vert_{L_{x}^2}^{1/2}\Vert\nabla\partial^{\alpha} b\Vert_{L_{x}^2}
     \\&\quad
     +
     C\Vert\nabla b\Vert_{L_{x}^{2}}\Vert\partial^{\alpha} b\Vert_{L_{x}^2}^{1/2}\Vert\nabla\partial^{\alpha} b\Vert_{L_{x}^2}^{1/2}\Vert\nabla\partial^{\alpha} u\Vert_{L_{x}^2}
\end{align*}
where we used Young's inequality in the last step.
Similarly, $I_2$ is estimated as
\begin{align*}
     &
     I_2
     \leq
     \frac{C}{\nu^3}\Vert\partial^{\alpha} u\Vert_{L_{x}^2}^2
     +
     \frac{C}{\nu}\Vert\partial^{\alpha} u\Vert_{L_{x}^2}
     +
     \frac{\nu}{8}\Vert\nabla\partial^{\alpha} u\Vert_{L_{x}^2}^2.
\end{align*}
By Cauchy-Schwarz inequality, we obtain,
\begin{align*}
     I_3 \leq \frac{g}{2}\Vert\partial^{\alpha} u\Vert_{L_{x}^2}^2 + \frac{g}{2}\Vert\partial^{\alpha} b\Vert_{L_{x}^2}^2.
\end{align*}
For the terms $I_4$ and $I_5$, 
we proceed similarly to the estimates of $I_1$.
Namely, we have
\begin{align*}
     I_4 + I_5
     &
     \leq
     C
     \left(\frac{C}{\nu\eta}+\frac{C}{\nu}+\frac{C}{\eta^3}+\frac{C}{\eta}+\frac{C}{\nu^3}\right)
     \left( \Vert\partial^{\alpha} b\Vert_{L_{x}^2}^2
     +
     \Vert\partial^{\alpha} u\Vert_{L_{x}^2}^2
     \right)
     \\&\qquad
     +
     \left(\frac{C}{\eta}+\frac{C}{\nu}\right) 
     \left( \Vert\partial^{\alpha} u\Vert_{L_{x}^2}
     +
     \Vert\partial^{\alpha} b\Vert_{L_{x}^2}\right)
     +
     \frac{\nu}{8}\Vert\nabla\partial^{\alpha} u\Vert_{L_{x}^2}^2
     +
     \frac{\eta}{8}\Vert\nabla\partial^{\alpha} b\Vert_{L_{x}^2}^2.
\end{align*}
Finally, the term $I_6$ is bounded as
\begin{align*}
     &
     I_6
     \leq
     \left(\frac{C}{\kappa^3}+\frac{C}{\kappa}\right)\Vert\partial^{\alpha} \theta\Vert_{L_{x}^2}^2
     +
     \frac{C}{\kappa}\Vert\partial^{\alpha} \theta\Vert_{L_{x}^2}
     +
     \frac{C}{\nu}\Vert\partial^{\alpha} u\Vert_{L_{x}^2}^2
     \\&\qquad
     +
     \frac{\nu}{8}\Vert\nabla\partial^{\alpha} u\Vert_{L_{x}^2}^2
     +
     \frac{\kappa}{2}\Vert\nabla\partial^{\alpha} \theta\Vert_{L_{x}^2}^2.
\end{align*}
Summing up the above estimates and denoting 
$$\bar{Q} = \Vert\partial^{\alpha} u\Vert_{L_{x}^2}^2 + \Vert\partial^{\alpha} b\Vert_{L_{x}^2}^2 + \Vert\partial^{\alpha} \theta\Vert_{L_{x}^2}^2,$$
we arrive at
\begin{equation} \label{H2_bound_Q_bar}
     \frac{d\bar{Q}}{dt}
     \leq
     C + C\bar{Q},
\end{equation}
where $C$ depends on $g$, $\nu$, $\eta$, $\kappa$, and $K_1(T')$ defined in \eqref{H1_bound_K1} (i.e., the bounds on the $H^1$ norms of $u$, $b$, and $\theta$).
Hence, by Gr\"onwall inequality,
we obtain $(u, b, \theta) \in L_{t}^{\infty}((0, T'); H^2\cap V)$. 
Also, we have 
$$\nu\int_{0}^{T'}\Vert\partial^{\alpha} u\Vert_{L_{x}^2}^2\,dt + \eta\int_{0}^{T'}\Vert\partial^{\alpha} b\Vert_{L_{x}^2}^2\,dt + \kappa\int_{0}^{T'}\Vert\partial^{\alpha} \theta\Vert_{L_{x}^2}^2\,dt$$
remains finite for $|\alpha|=2$. 
Next, we apply $\partial^{\alpha}$
with $|\alpha|=3$ to (\ref{S2}), 
and multiply the equations for $u$, $b$, 
and $\theta$ by $\partial^{\alpha} u$, $\partial^{\alpha} b$, and $\partial^{\alpha} \theta$, respectively, and get
\begin{align*}
  \left\{
     \begin{aligned}
     &
     \frac{1}{2}\frac{d}{dt}\Vert\partial^{\alpha} u\Vert_{L_{x}^2}^2
     +
     \nu\Vert\nabla\partial^{\alpha} u\Vert_{L_{x}^2}^2
     &&=
     \int_{\mathbb{T}^3}\partial^{\alpha}((b\cdot\nabla) b)\partial^{\alpha} u\,dx
     -
     \int_{\mathbb{T}^3}\partial^{\alpha}((u\cdot\nabla) u)\partial^{\alpha} u\,dx
     +
     g\int_{\mathbb{T}^3}\partial^{\alpha}\theta \partial^{\alpha} u\,dx
     \\&
     &&=
     J_1 + J_2 + J_3,
     \\&
     \frac{1}{2}\frac{d}{dt}\Vert\partial^{\alpha} b\Vert_{L_{x}^2}^2
     +
     \eta\Vert\nabla\partial^{\alpha} b\Vert_{L_{x}^2}^2
     &&=
     \int_{\mathbb{T}^3}\partial^{\alpha}((b\cdot\nabla) u)\partial^{\alpha} b\,dx
     -
     \int_{\mathbb{T}^3}\partial^{\alpha}((u\cdot\nabla) b)\partial^{\alpha} b\,dx
     =
     J_4 + J_5,
     \\&
     \frac{1}{2}\frac{d}{dt}\Vert\partial^{\alpha} \theta\Vert_{L_{x}^2}^2
     +
     \kappa\Vert\nabla\partial^{\alpha} \theta\Vert_{L_{x}^2}^2
     &&=
     -\int_{\mathbb{T}^3}\partial^{\alpha}((u\cdot\nabla) \theta)\partial^{\alpha} \theta\,dx
     =
     J_6.
   \end{aligned}
 \right.
\end{align*}
In order to estimate $J_1$, 
we apply Lemma~\ref{PL1} and obtain 
\begin{align*}
     J_1
     &
     \leq
     \sum_{0\leq|\zeta|\leq|\alpha|}\binom{\alpha}{\zeta}\int_{\mathbb{T}^3}|\partial^{\zeta}b| |\nabla\partial^{\alpha-\zeta}b| |\partial^{\alpha} u|\,dx
     \\&
     \leq
     C\Vert\nabla b\Vert_{L_{x}^2}\Vert\partial^{\alpha} u\Vert_{L_{x}^2}^{1/2}\Vert\nabla\partial^{\alpha} u\Vert_{L_{x}^2}^{1/2}\Vert\nabla\partial^{\alpha} b\Vert_{L_{x}^2}
     +
     C\sum_{|\zeta|=1}\Vert\partial^{\zeta} b\Vert_{L_{x}^2}^{1/2}\Vert\nabla\partial^{\zeta} b\Vert_{L_{x}^2}^{1/2}\Vert\partial^{\alpha} u\Vert_{L_{x}^2}\Vert\nabla\partial^{\alpha-\zeta} b\Vert_{L_{x}^2}
     \\&\qquad
     +
     C\sum_{|\zeta|=2}\Vert\partial^{\zeta} b\Vert_{L_{x}^2}^{3/2}\Vert\partial^{\alpha} b\Vert_{L_{x}^2}^{1/2}\Vert\nabla\partial^{\alpha} u\Vert_{L_{x}^2}
     +
     C\Vert\partial^{\alpha} b\Vert_{L_{x}^2}\Vert\partial^{\alpha} b\Vert_{L_{x}^2}^{1/2}\Vert\nabla\partial^{\alpha} b\Vert_{L_{x}^2}^{1/2}\Vert\partial^{\alpha} u\Vert_{L_{x}^2}
     \\&
     \leq
     \left(\frac{C}{\nu\eta}+\frac{C}{\eta}\right)\Vert\partial^{\alpha} u\Vert_{L_{x}^2}^2
     +
     \left(\frac{C}{\nu}+\frac{C}{\eta}\right)\Vert\partial^{\alpha} b\Vert_{L_{x}^2}
     +
     \frac{\nu}{8}\Vert\nabla\partial^{\alpha} u\Vert_{L_{x}^2}^2
     +
     \frac{\eta}{8}\Vert\nabla\partial^{\alpha} b\Vert_{L_{x}^2}^2,
\end{align*}
where we employed Young's inequality in the last inequality.
The estimates for $J_2$ are similar, i.e., we have
\begin{align*}
     &
     J_2
     \leq
     \frac{C}{\nu^3}\Vert\partial^{\alpha} u\Vert_{L_{x}^2}^2
     +
     \frac{C}{\nu}\Vert\partial^{\alpha} u\Vert_{L_{x}^2}
     +
     \frac{\nu}{8}\Vert\nabla\partial^{\alpha} u\Vert_{L_{x}^2}^2.
\end{align*}
Using Cauchy-Schwarz inequality, we obtain
\begin{align*}
     J_3\leq \frac{g}{2}\Vert\partial^{\alpha} u\Vert_{L_{x}^2}^2 + \frac{g}{2}\Vert\partial^{\alpha} b\Vert_{L_{x}^2}^2.
\end{align*}
Regarding $J_4$ and $J_5$, 
the estimates are similar to that of $J_1$.
Namely, we have
\begin{align*}
     J_4 + J_5
     &
     \leq
     C
     \left(\frac{C}{\nu\eta}+\frac{C}{\nu}+\frac{C}{\eta^3}+\frac{C}{\eta}+\frac{C}{\nu^3}\right)
     \left( \Vert\partial^{\alpha} b\Vert_{L_{x}^2}^2
     +
     \Vert\partial^{\alpha} u\Vert_{L_{x}^2}^2
     \right)
     \\&\qquad
     +
     \left(\frac{C}{\eta}+\frac{C}{\nu}\right) 
     \left( \Vert\partial^{\alpha} u\Vert_{L_{x}^2}
     +
     \Vert\partial^{\alpha} b\Vert_{L_{x}^2}\right)
     +
     \frac{\nu}{8}\Vert\nabla\partial^{\alpha} u\Vert_{L_{x}^2}^2
     +
     \frac{\eta}{8}\Vert\nabla\partial^{\alpha} b\Vert_{L_{x}^2}^2.
\end{align*}
Similarly, the term $J_6$ can be bounded as
\begin{align*}
     J_6
     &
     \leq
     \left(\frac{C}{\kappa^3}+\frac{C}{\kappa}\right)\Vert\partial^{\alpha} \theta\Vert_{L_{x}^2}^2
     +
     \frac{C}{\kappa}\Vert\partial^{\alpha} \theta\Vert_{L_{x}^2}
     +
     \frac{C}{\nu}\Vert\partial^{\alpha} u\Vert_{L_{x}^2}^2
     \\&\qquad
     +
     \frac{\nu}{8}\Vert\nabla\partial^{\alpha} u\Vert_{L_{x}^2}^2
     +
     \frac{\kappa}{2}\Vert\nabla\partial^{\alpha} \theta\Vert_{L_{x}^2}^2.
\end{align*}
Adding the above estimates and denoting 
$$Q = \Vert\partial^{\alpha} u\Vert_{L_{x}^2}^2 + \Vert\partial^{\alpha} b\Vert_{L_{x}^2}^2 + \Vert\partial^{\alpha} \theta\Vert_{L_{x}^2}^2,$$
we have
\begin{align*}
     \frac{dQ}{dt}
     \leq
     C + CQ,
\end{align*}
where $C$ depends on $g$, $\nu$, $\eta$, $\kappa$, and the bounds on the $H^2$ norms of $u$, $b$, and $\theta$.
Hence, using Gr\"onwall's inequality and combining all the above estimates, 
we finally obtain $(u, b, \theta) \in L_{t}^{\infty}((0, T'); H^3\cap V)$.
Furthermore, we have 
$$\nu\int_{0}^{T'}\Vert\nabla\partial^{\alpha} u\Vert_{L_{x}^2}^2\,dt + \eta\int_{0}^{T'}\Vert\nabla\partial^{\alpha} b\Vert_{L_{x}^2}^2\,dt + \kappa\int_{0}^{T'}\Vert\nabla\partial^{\alpha} \theta\Vert_{L_{x}^2}^2\,dt$$
remains finite for $|\alpha|=3$, i.e., $(u, b, \theta) \in L_{t}^2((0, T'); H^4\cap V)$.
Therefore, by slightly modifying the proof of the uniqueness of the non-diffusive system below, we obtain the uniqueness of the solution and Theorem~\ref{T2} is thus proven.
\endproof

\section{Uniqueness for the non-diffusive system}\label{sec4}
{\smallskip\noindent {\em Proof of uniqueness in Theorem~\ref{T3}.}}
In order to prove uniqueness, 
we use the fact that $(u, b, \theta) \in L^{\infty}([0, T^{*}); H^{m})$. 
Suppose that $(u^{(1)}, b^{(1)}, \theta^{(1)})$ and $(u^{(2)}, b^{(2)}, \theta^{(2)})$ 
are two solutions to the non-diffusive MHD-Boussinesq system (\ref{S1}).
By subtracting the two systems for the two solutions 
denoting $\tilde{u} = u^{(1)} - u^{(2)}$, $\tilde{p} = p^{(1)} - p^{(2)}$, $\tilde{b} = b^{(1)} - b^{(2)}$, and $\tilde{\theta} = \theta^{(1)} - \theta^{(2)}$, 
and by using H\"older's inequality, Gagliardo-Nirenberg-Sobolev inequality, and Young's inequality, 
to obtain
\begin{align*}
  \left\{
    \begin{aligned}
     &
     \frac{\partial\tilde{u}}{\partial t}
     -
     \nu\Delta\tilde{u}
     +
     (\tilde{u}\cdot\nabla)u^{(1)}
     +
     (u^{(2)}\cdot\nabla)\tilde{u}
     +
     \nabla\tilde{p}
     =
     (\tilde{b}\cdot\nabla)b^{(1)}
     +
     (b^{(2)}\cdot\nabla)\tilde{b}
     +
     g\tilde{\theta}e_{3},
     \\&
     \frac{\partial\tilde{b}}{\partial t}
     -
     \eta\Delta\tilde{b}
     +
     (\tilde{u}\cdot\nabla)b^{(1)}
     +
     (u^{(2)}\cdot\nabla)\tilde{b}
     =
     (\tilde{b}\cdot\nabla)u^{(1)}
     +
     (b^{(2)}\cdot\nabla)\tilde{u},
     \\&
     \frac{\partial\tilde{\theta}}{\partial t}
     +
     (\tilde{u}\cdot\nabla)\theta^{(1)}
     +
     (u^{(2)}\cdot\nabla)\tilde{\theta}
     =
     0,
   \end{aligned}
 \right.
\end{align*}
with $\nabla\cdot\tilde{u} = 0 = \nabla\tilde{b}$.
Multiply the above equations by $\tilde{u}$, $\tilde{b}$, and $\tilde{\theta}$, respectively, 
integrate over $\mathbb{T}^3$, and add, we get
\begin{align*}
     &\quad
     \frac{1}{2}\frac{d}{dt}\left(\Vert\tilde{u}\Vert_{L_{x}^2}^2 + \Vert\tilde{b}\Vert_{L_{x}^2}^2 + \Vert\tilde{\theta}\Vert_{L_{x}^2}^2\right)
     +
     \nu\Vert\nabla\tilde{u}\Vert_{L_{x}^2}^2
     +
     \eta\Vert\nabla\tilde{b}\Vert_{L_{x}^2}^2
     \\&
     =
     \int_{\mathbb{T}^3}(\tilde{u}\cdot\nabla)u^{(1)}\tilde{u}\,dx
     -
     \int_{\mathbb{T}^3}(\tilde{b}\cdot\nabla)b^{(1)}\tilde{u}\,dx
     +
     \int_{\mathbb{T}^3}g\tilde{\theta}e_{3}\tilde{u}\,dx
     \\&\qquad
     +
     \int_{\mathbb{T}^3}(\tilde{u}\cdot\nabla)b^{(1)}\tilde{b}\,dx
     -
     \int_{\mathbb{T}^3}(\tilde{b}\cdot\nabla)u^{(1)}\tilde{b}\,dx
     +
     \int_{\mathbb{T}^3}(\tilde{u}\cdot\nabla)\theta^{(1)}\tilde{\theta}\,dx
     \\&
     \leq
     C\Vert\nabla u^{(1)}\Vert_{L_{x}^2}\Vert\tilde{u}\Vert_{L_{x}^2}^{1/2}\Vert\nabla\tilde{u}\Vert_{L_{x}^2}^{3/2}
     +
     C\Vert\nabla b^{(1)}\Vert_{L_{x}^2}\Vert\tilde{b}\Vert_{L_{x}^2}^{1/2}\Vert\nabla\tilde{b}\Vert_{L_{x}^2}^{1/2}\Vert\nabla\tilde{u}\Vert_{L_{x}^2}
     +
     g\Vert\tilde{u}\Vert_{L_{x}^2}\Vert\tilde{\theta}\Vert_{L_{x}^2}
     \\&\qquad
     +
     C\Vert\nabla b^{(1)}\Vert_{L_{x}^2}\Vert\tilde{u}\Vert_{L_{x}^2}^{1/2}\Vert\nabla\tilde{u}\Vert_{L_{x}^2}^{1/2}\Vert\nabla\tilde{b}\Vert_{L_{x}^2}
     +
     C\Vert\nabla u^{(1)}\Vert_{L_{x}^2}\Vert\tilde{b}\Vert_{L_{x}^2}^{1/2}\Vert\nabla\tilde{b}\Vert_{L_{x}^2}^{3/2}
     \\&\qquad
     +
     C\Vert\tilde{u}\Vert_{L_{x}^2}^{1/2}\Vert\nabla\tilde{u}\Vert_{L_{x}^2}^{1/2}\Vert\nabla\nabla\theta^{(1)}\Vert_{L_{x}^2}\Vert\tilde{\theta}\Vert_{L_{x}^2}
     \\&
     \leq
     \frac{C}{\nu^3}\Vert\tilde{u}\Vert_{L_{x}^2}^2
     +
     \frac{\nu}{16}\Vert\nabla\tilde{u}\Vert_{L_{x}^2}^2
     +
     \frac{C}{\nu\eta}\Vert\tilde{b}\Vert_{L_{x}^2}^2
     +
     \frac{\nu}{16}\Vert\nabla\tilde{u}\Vert_{L_{x}^2}^2
     +
     \frac{\eta}{16}\Vert\nabla\tilde{b}\Vert_{L_{x}^2}^2
     \\&\qquad
     +
     \frac{g}{2}\Vert\tilde{\theta}\Vert_{L_{x}^2}^2
     +
     \frac{g}{2}\Vert\tilde{u}\Vert_{L_{x}^2}^2
     +
     \frac{C}{\nu\eta}\Vert\tilde{b}\Vert_{L_{x}^2}^2
     +
     \frac{\eta}{16}\Vert\nabla\tilde{b}\Vert_{L_{x}^2}^2
     +
     \frac{\nu}{16}\Vert\nabla\tilde{u}\Vert_{L_{x}^2}^2
     \\&\qquad
     +
     \frac{C}{\eta^3}\Vert\tilde{b}\Vert_{L_{x}^2}^2
     +
     \frac{\eta}{16}\Vert\nabla\tilde{b}\Vert_{L_{x}^2}^2
     +
     \frac{C}{\nu}\Vert\tilde{u}\Vert_{L_{x}^2}^2
     +
     \frac{\nu}{16}\Vert\nabla\tilde{u}\Vert_{L_{x}^2}^2
     +
     C\Vert\tilde{\theta}\Vert_{L_{x}^2}^2,
\end{align*}
where we used the bound in (\ref{H1_bound_K1}) and (\ref{H2_bound_Q_bar}) on $[0, T]$ for $T < T^{*}$. 
Let us denote $$X(t) = \Vert\tilde{u}\Vert_{L_{x}^2}^2 + \Vert\tilde{b}\Vert_{L_{x}^2}^2 + \Vert\tilde{\theta}\Vert_{L_{x}^2}^2,$$
for $0 \leq t \leq T < T^{*}$.
Then we have
\begin{align*}
     \frac{dX(t)}{dt}
     \leq
     CX(t),
\end{align*}
Gr\"onwall's inequality then gives continuity in the $L^\infty(0,T;L^2)$ norm.  Integrating, we also obtain continuity in the  $L^2(0,T;V)$ norm.  If the initial data is the same, then $X(0) = 0$, so we obtain uniqueness of the solutions
\endproof

\section{Proof of the regularity criterion}\label{sec5}
{\smallskip\noindent {\em Proof of Theorem~\ref{T0}.}}
We start by introducing the following notation. 
For the time interval $0 \leq t_1 < t_2 < \infty$, 
we denote
\begin{align*}
     (J(t_2))^2
     &:=
     \sup_{\tau \in (t_1, t_2)} \left\{ \Vert \nabla_{h} u(\tau) \Vert_{2}^2
     +
     \Vert \nabla_{h} b(\tau) \Vert_{2}^2
     \right\}
     +
     \int_{t_1}^{t_2} \Vert \nabla \nabla_{h} u(\tau)\Vert_{2}^2
     +
     \Vert \nabla \nabla_{h} b(\tau)\Vert_{2}^2\,d\tau,
\end{align*}
(recall that $\nabla_{h} = (\partial_1, \partial_2)$, and  $\Delta_{h} = \partial_{11} + \partial_{22}$).  We also denote
\begin{align*}
     (L(t_2))^2
     &:=
     \sup_{\tau \in (t_1, t_2)} \left\{ \Vert \partial_{3} u(\tau) \Vert_{2}^2
     +
     \Vert \partial_{3} b(\tau) \Vert_{2}^2
     \right\}
     +
     \int_{t_1}^{t_2} \Vert \nabla \partial_{3} u(\tau)\Vert_{2}^2
     +
     \Vert \nabla \partial_{3} b(\tau)\Vert_{2}^2\,d\tau.
\end{align*}
Aiming at a proof by contradiction, we denote the maximum time of existence and uniqueness of smooth solutions by $$T_{\text{max}} := \sup{\{ t \ge 0 | (u, b, \theta) \text{ is smooth on } (0, t)\}}.$$ 
Since $u_0$, $b_0$, and $\theta_0$ are in $H_{x}^3$, $T_{\text{max}}\in(0,\infty]$
If $T_{\text{max}} = \infty$, the proof is done.  Thus, we suppose $T_{\text{max}} < \infty$, and show that the solution can be extended beyond $T_{\text{max}}$, which is a contradiction.  
First, we choose $\epsilon > 0$ sufficiently small, 
say, $\epsilon < 1/(16C_{max})$, 
where $C_{max}$ is the maximum of all the constants in the following argument, 
depending on the space dimension, the constant $g$, 
the first eigenvalue $\lambda_1$ of the operator $-\Delta$, as well as the spatial-temporal $L^2$-norm of the solution up to $T_{max}$. 
Then, we fix $T_1 \in (0, T_{\text{max}})$ such that 
$T_{\text{max}} - T_1 < \epsilon$, 
and
\begin{equation}
     \int_{T_1}^{T_{\text{max}}} \Vert \nabla u(\tau) \Vert_{L_{x}^2}^2 + \Vert \nabla b(\tau) \Vert_{L_{x}^2}^2 + \Vert\theta\Vert_{L_{x}^2}^2\,d\tau < \epsilon,
     \label{Ineq}
\end{equation}
as well as 
\begin{equation}
     \int_{T_1}^{T_{\text{max}}}\Vert u_2(\tau) \Vert_{L_{x}^{s}}^{r} + \Vert u_3(\tau) \Vert_{L_{x}^{s}}^{r} + \Vert b_2(\tau) \Vert_{L_{x}^{s}}^{r} + \Vert b_3(\tau) \Vert_{L_{x}^{s}}^{r}\,d\tau < \epsilon.
\end{equation}
We see that the proof is complete if we show that $\Vert \nabla u(T_2) \Vert_{2}^2 + \Vert \nabla b(T_2) \Vert_{2}^2 + \Vert \nabla \theta(T_2) \Vert_{2}^2 \leq C < \infty$, 
for any $T_2 \in (T_1, T_{\text{max}})$ and $C$ in independent of the choice of $T_2$.
In fact, due to the continuity of integral, we can extend the the regularity of $u$ beyond $T_{\text{max}}$ 
and this becomes a contradiction to the definition of $T_{\text{max}}$.
Therefore, it is sufficient to prove that 
$J(T_2)^2 + L(T_2)^2 \leq C < \infty$ in view of the equation for $\theta$ in (\ref{S1}) 
for some constant $C$ independent of $T_2$.
We take the approach of \cite{S}, which first bounds $L(T_2)$ by $J(T_2)$, 
then closes the estimates by obtaining an uniform upper bound on the latter. 
The regularity of $\theta$ thus follows from the higher order regularity of $u$ and $b$.
To start, we multiply the equations for $u$ and $b$ 
in (\ref{S1}) by $-\partial_{33}^{2}u$ and $-\partial_{33}^{2}b$ respectively, 
integrate over $\mathbb{T}^3 \times (T_1, T_2)$, and sum to obtain
\begin{align*}
     &\quad
     \frac{1}{2}\left(\Vert \partial_{3}u(T_2)\Vert_{L_{x}^2}^2
     +
     \Vert \partial_{3}b(T_2)\Vert_{L_{x}^2}^2
     \right)
     \\&\qquad
     +
     \int_{T_1}^{T_2} \int_{\mathbb{T}^3} \nu\Vert \nabla\partial_{3}u\Vert_{L_{x}^2}^2
     +
     \eta\Vert \nabla\partial_{3}b\Vert_{L_{x}^2}^2\,dx\,d\tau
     \\&
     =
     \frac{1}{2}\left(\Vert \partial_{3}u(T_1)\Vert_{L_{x}^2}^2
     +
     \Vert \partial_{3}b(T_1)\Vert_{L_{x}^2}^2 
     \right)
     \\&\qquad
     -
     \sum_{j, k = 1}^{3}\int_{T_{1}}^{T_{2}}\int_{\mathbb{T}^3} \partial_{3}u_{j}\partial_{j}u_{k}\partial_{3}u_{k}\,dx\,d\tau
     +
     \sum_{j, k = 1}^{3}\int_{T_{1}}^{T_{2}}\int_{\mathbb{T}^3} \partial_{3}b_{j}\partial_{j}b_{k}\partial_{3}u_{k}\,dx\,d\tau
     \\&\qquad
     -
     \sum_{j, k = 1}^{3}\int_{T_{1}}^{T_{2}}\int_{\mathbb{T}^3} \partial_{3}u_{j}\partial_{j}b_{k}\partial_{3}b_{k}\,dx\,d\tau
     +
     \sum_{j, k = 1}^{3}\int_{T_{1}}^{T_{2}}\int_{\mathbb{T}^3} \partial_{3}b_{j}\partial_{j}u_{k}\partial_{3}b_{k}\,dx\,d\tau
     \\&\qquad
     -
     g\sum_{k = 1}^{3}\int_{T_{1}}^{T_{2}}\int_{\mathbb{T}^3} \theta e_3 \partial_{33}u_{k}\,dx\,d\tau,
\end{align*}
where we used the divergence-free condition and Lemma~\ref{L02}. Then we denote the last five integrals on the right side of the above equation by 
$\uppercase\expandafter{\romannumeral1}$, $\uppercase\expandafter{\romannumeral2}$, $\uppercase\expandafter{\romannumeral3}$, $\uppercase\expandafter{\romannumeral4}$, and $\uppercase\expandafter{\romannumeral5}$, respectively.
In order to estimate $\uppercase\expandafter{\romannumeral1}$ we first rewrite it as
\begin{align*}
     \uppercase\expandafter{\romannumeral1}
     &
     =
     -\sum_{j, k = 1}^{2}\int_{T_{1}}^{T_{2}}\int_{\mathbb{T}^3} \partial_{3}u_{j}\partial_{j}u_{k}\partial_{3}u_{k}\,dx\,d\tau
     -
     \sum_{j = 1}^{2}\int_{T_{1}}^{T_{2}}\int_{\mathbb{T}^3} \partial_{3}u_{j}\partial_{j}u_{3}\partial_{3}u_{3}\,dx\,d\tau
     \\&\qquad
     -
     \sum_{k = 1}^{2}\int_{T_{1}}^{T_{2}}\int_{\mathbb{T}^3} \partial_{3}u_{3}\partial_{3}u_{k}\partial_{3}u_{k}\,dx\,d\tau
     -
     \int_{T_{1}}^{T_{2}}\int_{\mathbb{T}^3} \partial_{3}u_{3}\partial_{3}u_{3}\partial_{3}u_{3}\,dx\,d\tau
     \\&
     =
     \sum_{j, k = 1}^{2}\int_{T_{1}}^{T_{2}}\int_{\mathbb{T}^3} u_{k}\left(\partial_{3}u_{k}\partial_{3j}^{2}u_{j}
     +
     \partial_{3}u_{j}\partial_{3j}^{2}u_{k}\right)\,dx\,d\tau
     -
     \uppercase\expandafter{\romannumeral1}_{a}
     -
     \uppercase\expandafter{\romannumeral1}_{b}
     -
     \uppercase\expandafter{\romannumeral1}_{c}.
\end{align*}
By Lemma~\ref{PL1}, 
the first two integrals on the right side of $\uppercase\expandafter{\romannumeral1}$ are bounded by
\begin{align*}
     &\quad
     C\int_{T_{1}}^{T_{2}}\int_{\mathbb{T}^3} |u| |\partial_{3}u| |\nabla\partial_{3}u|\,dx\,d\tau
     \\&
     \leq
     C\int_{T_{1}}^{T_{2}}\Vert u \Vert_{L_{x}^{6}}\Vert\partial_{3}u\Vert_{L_{x}^3}\Vert\nabla_{h}\partial_{3}u\Vert_{L_{x}^2}\,d\tau
     \\&
     \leq
     C\int_{T_{1}}^{T_{2}}\Vert u \Vert_{L_{x}^{6}}\Vert\partial_{3}u\Vert_{L_{x}^2}^{\frac{1}{2}}\Vert\partial_{3}u\Vert_{L_{x}^6}^{\frac{1}{2}}\Vert\nabla_{h}\partial_{3}u\Vert_{L_{x}^2}\,d\tau
     \\&
     \leq
     C\Vert \nabla_{h} u\Vert_{L_{t}^{\infty}L_{x}^{2}}^{\frac{2}{3}}\Vert \partial_{3} u\Vert_{L_{t}^{\infty}L_{x}^{2}}^{\frac{1}{3}}\Vert\partial_{3}u\Vert_{L_{t}^{2}L_{x}^2}^{\frac{1}{2}}\Vert\nabla_{h}\partial_{3}u\Vert_{L_{t}^{2}L_{x}^2}^{\frac{1}{3}}\Vert\partial_{33}^2u\Vert_{L_{t}^{2}L_{x}^2}^{\frac{1}{6}}\Vert\nabla_{h}\partial_{3}u\Vert_{L_{t}^2L_{x}^2}
     \\&
     \leq
     C\epsilon L^{\frac{1}{2}}(T_2) J^2(T_2),
\end{align*}
where the $L^\infty_t$ norms are taken over the interval $(T_1,T_2)$ 
and we used Lemma~\ref{PL2} in the second to the last inequality. 
Regarding $\uppercase\expandafter{\romannumeral1}_{a}$, $\uppercase\expandafter{\romannumeral1}_{b}$, and $\uppercase\expandafter{\romannumeral1}_{c}$, 
we first integrate by parts, then estimate as
\begin{align*}
     \uppercase\expandafter{\romannumeral1}_{a}
     +
     \uppercase\expandafter{\romannumeral1}_{b}
     +
     \uppercase\expandafter{\romannumeral1}_{c}
     &
     =
     \sum_{j= 1}^{2}\int_{T_{1}}^{T_{2}}\int_{\mathbb{T}^3} u_{3}\partial_{3}u_{j}\partial_{3j}^{2}u_{3}\,dx\,d\tau
     +
     \sum_{j= 1}^{2}\int_{T_{1}}^{T_{2}}\int_{\mathbb{T}^3} u_{3}\partial_{j}u_{3}\partial_{33}^{2}u_{j}\,dx\,d\tau
     \\&\qquad
     +
     2\sum_{k = 1}^{2}\int_{T_{1}}^{T_{2}}\int_{\mathbb{T}^3} u_{3}\partial_{3}u_{k}\partial_{33}u_{k}\,dx\,d\tau
     +
     2\int_{T_{1}}^{T_{2}}\int_{\mathbb{T}^3} u_{3}\partial_{3}u_{3}\partial_{33}u_{3}\,dx\,d\tau
     \\&
     \leq
     C\int_{T_{1}}^{T_{2}} |u_{3}| |\nabla_{h}u| |\nabla\partial_{3}u|\,dx\,d\tau
     +
     C\int_{T_{1}}^{T_{2}} |u_{3}| |\partial_{3}u| |\nabla\partial_{3}u|\,dx\,d\tau
     \\&
     \leq
     C\int_{T_{1}}^{T_{2}}\Vert u_{3}\Vert_{L_{x}^{s}}\Vert\nabla_{h}u\Vert_{L_{x}^2}^{1-\frac{3}{s}}\Vert\nabla\partial_{3}u\Vert_{L_{x}^2}^{1+\frac{3}{s}}\,d\tau
     +
     C\int_{T_{1}}^{T_{2}}\Vert u_{3}\Vert_{L_{x}^{s}}\Vert\partial_{3}u\Vert_{L_{x}^2}^{1-\frac{3}{s}}\Vert\nabla\partial_{3}u\Vert_{L_{x}^2}^{1+\frac{3}{s}}\,d\tau
     \\&
     \leq
     C(T_2 - T_1)^{1-(\frac{2}{r} + \frac{3}{s})}\Vert u_{3}\Vert_{L_{t}^{r}L_{x}^{s}}\Vert\nabla_{h}u\Vert_{L_{t}^{\infty}L_{x}^2}^{1-\frac{3}{s}}\Vert\nabla\partial_{3}u\Vert_{L_{t}^{2}L_{x}^2}^{1+\frac{3}{s}}
     \\&\qquad
     +
     C(T_2 - T_1)^{1-(\frac{2}{r} + \frac{3}{s})}\Vert u_{3}\Vert_{L_{t}^{r}L_{x}^{s}}\Vert\partial_{3}u\Vert_{L_{t}^{\infty}L_{x}^2}^{1-\frac{3}{s}}\Vert\nabla\partial_{3}u\Vert_{L_{t}^{2}L_{x}^2}^{1+\frac{3}{s}}
     \\&
     \leq
     C\epsilon J^{1-\frac{3}{s}}(T_2) L^{1+\frac{3}{s}}(T_2)
     +
     C\epsilon L^2(T_2),
\end{align*}
where we used the fact that $\Vert\nabla u\Vert_{L_{t}^2L_{x}^2}^{1/2}$ is small over the interval $(T_1, T_2)$ 
and the constant $C$ is independent of $T_2$.
Next, we estimate $\uppercase\expandafter{\romannumeral2}$. 
Proceeding similarly as the estimates for $\uppercase\expandafter{\romannumeral1}$, 
we first integrate by parts and rewrite $\uppercase\expandafter{\romannumeral2}$ as
\begin{align*}
     \uppercase\expandafter{\romannumeral2}
     &
     =
     \sum_{j = 1}^{3}\sum_{k = 1}^{2}\int_{T_{1}}^{T_{2}}\int_{\mathbb{T}^3} b_{k}\partial_{3}b_{j}\partial_{3j}^2u_{k}\,dx\,d\tau
     +
     \sum_{j = 1}^{3}\int_{T_{1}}^{T_{2}}\int_{\mathbb{T}^3} b_{3}\partial_{3}b_{j}\partial_{3j}^2u_{3}\,dx\,d\tau
     \\&
     \leq
     C\int_{T_{1}}^{T_{2}}\int_{\mathbb{T}^3} |b| |\partial_{3}b| |\nabla_{h}\partial_{3}u|\,dx\,d\tau
     +
     C\int_{T_{1}}^{T_{2}}\int_{\mathbb{T}^3} |b_3| |\partial_{3}b| |\nabla\partial_{3}u|\,dx\,d\tau. 
\end{align*}
Therefore, by Lemma~\ref{PL1} and Lemma~\ref{PL2}, we get
\begin{align*}
     \uppercase\expandafter{\romannumeral2}
     &
     \leq
     C\int_{T_{1}}^{T_{2}}\Vert b \Vert_{L_{x}^{6}}\Vert\partial_{3}b\Vert_{L_{x}^3}\Vert\nabla_{h}\partial_{3}u\Vert_{L_{x}^2}\,d\tau
     +
     C\int_{T_{1}}^{T_{2}}\int_{\mathbb{T}^3} (|u_3| + |b_3|) (|\partial_{3}u| + |\partial_{3}b|) (|\nabla\partial_{3}u| + |\nabla\partial_{3}b|)\,dx\,d\tau
     \\&
     \leq
     C\int_{T_{1}}^{T_{2}}\Vert b \Vert_{L_{x}^{6}}\Vert\partial_{3}b\Vert_{L_{x}^2}^{\frac{1}{2}}\Vert\nabla\partial_{3}b\Vert_{L_{x}^2}^{\frac{1}{2}}\Vert\nabla_{h}\partial_{3}u\Vert_{L_{x}^2}\,d\tau
     \\&\qquad
     +
     C\int_{T_{1}}^{T_{2}}(\Vert u_{3}\Vert_{L_{x}^{s}} + \Vert b_{3}\Vert_{L_{x}^{s}})(\Vert\partial_{3}u\Vert_{L_{x}^2} + \Vert\partial_{3}b\Vert_{L_{x}^2})^{1-\frac{3}{s}}(\Vert\nabla\partial_{3}u\Vert_{L_{x}^2} + \Vert\nabla\partial_{3}b\Vert_{L_{x}^2})^{1+\frac{3}{s}}\,d\tau
     \\&
     \leq
     C\Vert \nabla_{h} b\Vert_{L_{t}^{\infty}L_{x}^{2}}^{\frac{2}{3}}\Vert \partial_{3} b\Vert_{L_{t}^{\infty}L_{x}^{2}}^{\frac{1}{3}}\Vert\partial_{3}b\Vert_{L_{t}^{2}L_{x}^2}^{\frac{1}{2}}\Vert\nabla_{h}\partial_{3}b\Vert_{L_{t}^{2}L_{x}^2}^{\frac{1}{3}}\Vert\partial_{33}^2b\Vert_{L_{t}^{2}L_{x}^2}^{\frac{1}{6}}\Vert\nabla_{h}\partial_{3}u\Vert_{L_{t}^2L_{x}^2}
     \\&\qquad
     +
     C(T_2 - T_1)^{1-(\frac{2}{r} + \frac{3}{s})}(\Vert u_{3}\Vert_{L_{t}^{r}L_{x}^{s}} + \Vert b_{3}\Vert_{L_{t}^{r}L_{x}^{s}})(\Vert\partial_{3}u\Vert_{L_{t}^{\infty}L_{x}^2} +\Vert\partial_{3}b\Vert_{L_{t}^{\infty}L_{x}^2})^{1-\frac{3}{s}}
     \\&\quad\qquad
     \times
     (\Vert\nabla\partial_{3}u\Vert_{L_{t}^{2}L_{x}^2} + \Vert\nabla\partial_{3}b\Vert_{L_{t}^{2}L_{x}^2})^{1+\frac{3}{s}}
     \\&
     \leq
     C\epsilon L^{\frac{1}{2}}(T_2) J^2(T_2)
     +
     C\epsilon J^{1-\frac{3}{s}}(T_2) L^{1+\frac{3}{s}}(T_2)
     +
     C\epsilon L^2(T_2).
\end{align*}
The terms $\uppercase\expandafter{\romannumeral3}$ and $\uppercase\expandafter{\romannumeral4}$ are estimated analogously, 
i.e., we have
\begin{align*}
     &
     \uppercase\expandafter{\romannumeral3}
     +
     \uppercase\expandafter{\romannumeral4}
     \leq
     C\epsilon L^{\frac{1}{2}}(T_2) J^2(T_2)
     +
     C\epsilon J^{1-\frac{3}{s}}(T_2) L^{1+\frac{3}{s}}(T_2)
     +
     C\epsilon L^2(T_2),
\end{align*}
where the constant $C$ does not depend on $T_2$.
We estimate the term $\uppercase\expandafter{\romannumeral5}$ as
\begin{align*}
     \uppercase\expandafter{\romannumeral5}
     &
     =
     -\sum_{k=1}^{3}\int_{T_{1}}^{T_{2}}\int_{\mathbb{T}^3} \theta e_3\partial_{33}u_{k}\,d\tau
     \leq
     C\Vert\theta\Vert_{L_{x, t}^{2}} \Vert\partial_{33}u\Vert_{L_{x, t}^{2}}
     \leq
     C\Vert\theta_0\Vert_{L_{x}^{2}} \Vert\partial_{33}u\Vert_{L_{x, t}^{2}}
     \leq
     C\epsilon L(T_2).
\end{align*}
Collecting the above estimate for $\uppercase\expandafter{\romannumeral1}$ through $\uppercase\expandafter{\romannumeral5}$ 
and using Young's inequality, we obtain
\begin{align*}
     L^{2}(T_2)
     &
     \leq
     C
     +
     C\epsilon L^{\frac{1}{2}}(T_2)J^2(T_2)
     +
     C\epsilon L^{1+\frac{3}{s}}(T_2) J^{1-\frac{3}{s}}(T_2)
     +
     C\epsilon L^2(T_2)
     +
     C\epsilon L(T_2)
     \\&
     \leq
     C
     +
     C\epsilon L^{2}(T_2)
     +
     C\epsilon J^{\frac{8}{3}}(T_2)
     +
     C\epsilon J^2(T_2)
     +
     C\epsilon L(T_2).
\end{align*}
Thus, with our choice of $\epsilon >0$ earlier, we get
\begin{equation}
     L(T_2)
     \leq
     C
     +
     CJ(T_{2})^{\frac{4}{3}}.
     \label{Ineq01}
\end{equation}
Next, in order to bound $J(T_2)$, 
we multiply the equation for $u$ and $b$ in \ref{S1} 
by $-\Delta_{h}u$ and  
$-\Delta_{h}b$, respectively,
integrate over $\mathbb{T}^3 \times (T_1, T_2)$, sum up, 
integrate by parts and get
\begin{align*}
     &\quad
     \frac{1}{2}\left(\Vert \nabla_{h}u(T_2)\Vert_{L_{x}^2}^2
     +
     \Vert \nabla_{h}b(T_2)\Vert_{L_{x}^2}^2
     \right)
     +
     \int_{T_1}^{T_2} \int_{\mathbb{T}^3} \Vert \nabla\nabla_{h}u\Vert_{L_{x}^2}^2
     +
     \Vert \nabla\nabla_{h}b\Vert_{L_{x}^2}^2
     \\&
     =
     \frac{1}{2}\left(\Vert \nabla_{h}u(T_1)\Vert_{L_{x}^2}^2
     +
     \Vert \nabla_{h}b(T_1)\Vert_{L_{x}^2}^2\right)
     \\&\qquad
     -
     \sum_{j, k = 1}^{3}\sum_{i = 1}^{2}\int_{T_1}^{T_2}\int_{\mathbb{T}^3}\partial_{i}u_{j}\partial_{j}u_{k}\partial_{i}u_{k}\,dx\,d\tau
     +
     \sum_{j, k = 1}^{3}\sum_{i = 1}^{2}\int_{T_1}^{T_2}\int_{\mathbb{T}^3} \partial_{i}b_{j}\partial_{j}b_{k}\partial_{i}u_{k}\,dx\,d\tau
     \\&\qquad
     -
     \sum_{j, k = 1}^{3}\sum_{i = 1}^{2}\int_{T_1}^{T_2}\int_{\mathbb{T}^3} \partial_{i}u_{j}\partial_{j}b_{k}\partial_{i}b_{k}\,dx\,d\tau
     +
     \sum_{j, k = 1}^{3}\sum_{i = 1}^{2}\int_{T_1}^{T_2}\int_{\mathbb{T}^3} \partial_{i}b_{j}\partial_{j}u_{k}\partial_{i}b_{k}\,dx\,d\tau
     \\&\qquad
     -
     g\sum_{k = 1}^{3}\sum_{i = 1}^{2}\int_{T_1}^{T_2}\int_{\mathbb{T}^3} \theta e_3\partial_{ii}u_{k}\,dx\,d\tau,
\end{align*}
where we used the divergence-free condition and Lemma~\ref{L02}.
Denote by $\tilde{\uppercase\expandafter{\romannumeral1}}$ through $\tilde{\uppercase\expandafter{\romannumeral5}}$ the last five integrals on the right side of the above equation, respectively. Integrating by parts, we first rewrite $\tilde{\uppercase\expandafter{\romannumeral1}}$ as
\begin{align*}
     \tilde{\uppercase\expandafter{\romannumeral1}}
     &
     =
     -
     \sum_{i, j, k = 1}^{2}\int_{T_1}^{T_2}\int_{\mathbb{T}^3}\partial_{i}u_{j}\partial_{j}u_{k}\partial_{i}u_{k}\,dx\,d\tau
     -
     \sum_{i, j = 1}^{2}\int_{T_1}^{T_2}\int_{\mathbb{T}^3}\partial_{i}u_{j}\partial_{j}u_{3}\partial_{i}u_{3}\,dx\,d\tau
     \\&\qquad
     -
     \sum_{i, k = 1}^{2}\int_{T_1}^{T_2}\int_{\mathbb{T}^3}\partial_{i}u_{3}\partial_{3}u_{k}\partial_{i}u_{k}\,dx\,d\tau
     -
     \sum_{i = 1}^{2}\int_{T_1}^{T_2}\int_{\mathbb{T}^3}\partial_{i}u_{3}\partial_{3}u_{3}\partial_{i}u_{3}\,dx\,d\tau   
     \\&
     =
     \frac{1}{2}\sum_{j, k=1}^2 \int_{T_1}^{T_2}\int_{\mathbb{T}^3} u_{3}\partial_{j}u_{k}\partial_{3j}^{2}u_{k}\,dx\,d\tau
     -
     \int_{T_1}^{T_2}\int_{\mathbb{T}^3} u_{3}\partial_{1}u_{1}\partial_{32}^{2}u_{2}\,dx\,d\tau
     -
     \int_{T_1}^{T_2}\int_{\mathbb{T}^3} u_{3}\partial_{2}u_{2}\partial_{31}^{2}u_{1}\,dx\,d\tau
     \\&\qquad
     +
     \int_{T_1}^{T_2}\int_{\mathbb{T}^3} u_{3}\partial_{1}u_{2}\partial_{32}^{2}u_{1}\,dx\,d\tau
     +
      \int_{T_1}^{T_2}\int_{\mathbb{T}^3} u_{3}\partial_{2}u_{1}\partial_{31}^{2}u_{2}\,dx\,d\tau
     \\&\qquad
     +
     \sum_{i, j = 1}^{2}\int_{T_1}^{T_2}\int_{\mathbb{T}^3}u_{3}\partial_{i}u_{j}\partial_{3j}^{2}u_{3}\,dx\,d\tau
     +
     \sum_{i, j = 1}^{2}\int_{T_1}^{T_2}\int_{\mathbb{T}^3}u_{3}\partial_{j}u_{3}\partial_{3i}^{2}u_{j}\,dx\,d\tau
     \\&\qquad
     +
     \sum_{i, k = 1}^{2}\int_{T_1}^{T_2}\int_{\mathbb{T}^3}u_{3}\partial_{3}u_{k}\partial_{ii}^{2}u_{k}\,dx\,d\tau
     +
     \sum_{i, k = 1}^{2}\int_{T_1}^{T_2}\int_{\mathbb{T}^3}u_{3}\partial_{i}u_{k}\partial_{33}^{2}u_{k}\,dx\,d\tau
     \\&\qquad
     +
     2\sum_{i = 1}^{2}\int_{T_1}^{T_2}\int_{\mathbb{T}^3}u_{3}\partial_{i}u_{3}\partial_{3i}^{2}u_{3}\,dx\,d\tau, 
\end{align*} 
where we applied Lemma~\ref{L01} to the first term on the right side of the first equality above.
Thus, by H\"older and Sobolev inequalities, we bound $\tilde{\uppercase\expandafter{\romannumeral1}}$ as
\begin{align*}
     \tilde{\uppercase\expandafter{\romannumeral1}}
     &
     \leq
     C\int_{T_1}^{T_2}\int_{\mathbb{T}^3} |u_{3}| (|\nabla_{h}u| + |\partial_{3}u|)|\nabla\nabla_{h}u|\,dx\,d\tau
     \\&
     \leq
     C\int_{T_1}^{T_2} \Vert u_{3}\Vert_{L_{x}^{s}} \Vert\nabla_{h} u\Vert_{L_{x}^2}^{1-\frac{3}{s}} \Vert\nabla\nabla_{h}u\Vert_{L_{x}^2}^{1+\frac{3}{s}}\,d\tau
     \\&\qquad
     +
     C\int_{T_1}^{T_2}\Vert u_{3}\Vert_{L_{x}^{s}} \Vert\partial_{3} u\Vert_{L_{x}^2}^{1-\frac{3}{s}} \Vert\nabla_{h}\partial_{3}u\Vert_{L{x}^2}^{\frac{2}{s}} \Vert\partial_{33}^{2}u\Vert_{L{x}^2}^{\frac{1}{s}} \Vert\nabla\nabla_{h}u\Vert_{L_{x}^2}\,d\tau
     \\&
     \leq
     C(T_2 - T_1)^{1-(\frac{2}{r} + \frac{3}{s})} \Vert u_{3}\Vert_{L_{t}^{r}L_{x}^{s}} \Vert\nabla_{h} u\Vert_{L_{t}^{\infty}L_{x}^2}^{1-\frac{3}{s}} \Vert\nabla\nabla_{h}u\Vert_{L_{t}^2L_{x}^2}^{1+\frac{3}{s}}
     \\&\qquad
     +
     C(T_2 - T_1)^{1-(\frac{2}{r} + \frac{3}{s})} \Vert u_{3}\Vert_{L_{t}^{r}L_{x}^{s}} \Vert\partial_{3} u\Vert_{L_{t}^{2}L_{x}^2}^{\frac{s-2}{4s}}\Vert\partial_{3} u\Vert_{L_{t}^{\infty}L_{x}^2}^{\frac{3s-10}{4s}}\times \Vert\nabla\partial_{3}u\Vert_{L_{t}^2L_{x}^2}^{\frac{1}{s}}\Vert\nabla\nabla_{h}u\Vert_{L_{t}^2L_{x}^2}^{1+\frac{2}{s}}
     \\&
     \leq
     C
     +
     C\epsilon J^2(T_2)
     +
     CC\epsilon J^{\frac{4}{3}\frac{3s-6}{4s} + 1 + \frac{2}{s}}
     \\&
     \leq
     C
     +
     C\epsilon J^2(T_2),
\end{align*}
where we used (\ref{Ineq01}) and the fact that $T_2 - T_1 < \epsilon$ and $2/r + 3/s = 3/4 + 1/(2s)$ for $s > 10/3$.
In order to estimate $\tilde{\uppercase\expandafter{\romannumeral2}}$, 
we proceed a bit differently since Lemma~\ref{L01} is not available for convective terms mixed with $u$ and $b$.
Instead, we integrate by parts and use the divergence-free condition $\partial_{1}b_1 = -\partial_{2}b_2 - \partial_{3}b_3$ and obtain
\begin{align*}
     \tilde{\uppercase\expandafter{\romannumeral2}}
     &
     =
     \sum_{j, k = 1}^{3}\sum_{i = 1}^{2}\int_{T_1}^{T_2}\int_{\mathbb{T}^3} \partial_{i}b_{j}\partial_{j}b_{k}\partial_{i}u_{k}\,dx\,d\tau
     \\&
     =
     \sum_{i = 1}^{2} \int_{T_1}^{T_2}\int_{\mathbb{T}^3} \partial_{j}b_{1}\partial_{1}b_{1}\partial_{i}u_{1}\,dx\,d\tau
     +
     \sum_{i = 1}^{2}\sum_{k = 2}^{3}\int_{T_1}^{T_2}\int_{\mathbb{T}^3} \partial_{i}b_{1}\partial_{1}b_{k}\partial_{i}u_{k}\,dx\,d\tau
     \\&\qquad
     +
     \sum_{i = 1}^{2}\sum_{k = 1}^{3}\sum_{j = 2}^{3}\int_{T_1}^{T_2}\int_{\mathbb{T}^3} \partial_{i}b_{j}\partial_{j}b_{k}\partial_{i}u_{k}\,dx\,d\tau
     \\&
     =
     \sum_{i = 1}^{2} \int_{T_1}^{T_2}\int_{\mathbb{T}^3} \partial_{j}b_{1}(-b_{2}\partial_{2} - b_{3}\partial_{3})\partial_{i}u_{1}\,dx\,d\tau
     \\&\qquad
     -
     \sum_{i = 1}^{2}\sum_{k = 2}^{3}\int_{T_1}^{T_2}\int_{\mathbb{T}^3} u_{k}\partial_{i}b_{1}\partial_{1i}^{2}b_{k}\,dx\,d\tau
     -
     \sum_{i = 1}^{2}\sum_{k = 2}^{3}\int_{T_1}^{T_2}\int_{\mathbb{T}^3} u_{k}\partial_{1}b_{k}\partial_{ii}^{2}b_{1}\,dx\,d\tau
     \\&\qquad
     -
     \sum_{i = 1}^{2}\sum_{k = 1}^{3}\sum_{j = 2}^{3}\int_{T_1}^{T_2}\int_{\mathbb{T}^3} b_{j}\partial_{j}b_{k}\partial_{ii}^{2}u_{k}\,dx\,d\tau
     -
     \sum_{i = 1}^{2}\sum_{k = 1}^{3}\sum_{j = 2}^{3}\int_{T_1}^{T_2}\int_{\mathbb{T}^3} b_{j}\partial_{i}u_{k}\partial_{ij}^{2}b_{k}\,dx\,d\tau.
\end{align*}
Then after integration by parts to the first term on the right side of the above equation,
we bound $\tilde{\uppercase\expandafter{\romannumeral2}}$ as
\begin{align*}
     \tilde{\uppercase\expandafter{\romannumeral2}}
     &
     \leq
     C\int_{T_1}^{T_2}\int_{\mathbb{T}^3} (|b_{2}| + |b_{3}|) (|\nabla_{h}u| + |\nabla_{h}b| + |\partial_{3}u| + |\partial_{3}b|)(|\nabla\nabla_{h}u| + |\nabla\nabla_{h}b|)\,dx\,d\tau
     \\&
     \leq
     C\int_{T_1}^{T_2} (\Vert b_{2}\Vert_{L_{x}^{s}} + \Vert b_{3}\Vert_{L_{x}^{s}}) (\Vert\nabla_{h} u\Vert_{L_{x}^2} + \Vert\nabla_{h} b\Vert_{L_{x}^2})^{1-\frac{3}{s}} (\Vert\nabla\nabla_{h}u\Vert_{L_{x}^2} + \Vert\nabla\nabla_{h}b\Vert_{L_{x}^2})^{1+\frac{3}{s}}\,d\tau
     \\&\qquad
     +
     C\int_{T_1}^{T_2}(\Vert b_{2}\Vert_{L_{x}^{s}} + \Vert b_{3}\Vert_{L_{x}^{s}}) (\Vert\partial_{3} u\Vert_{L_{x}^2} + \Vert\partial_{3} b\Vert_{L_{x}^2})^{1-\frac{3}{s}} (\Vert\nabla_{h}\partial_{3}u\Vert_{L{x}^2} + \Vert\nabla_{h}\partial_{3}b\Vert_{L{x}^2})^{\frac{2}{s}} 
     \\&\quad\qquad
     \times
     (\Vert\partial_{33}^{2}u\Vert_{L{x}^2} + \Vert\partial_{33}^{2}b\Vert_{L{x}^2})^{\frac{1}{s}} (\Vert\nabla\nabla_{h}u\Vert_{L_{x}^2} + \Vert\nabla\nabla_{h}b\Vert_{L_{x}^2})\,d\tau
     \\&
     \leq
     C(T_2 - T_1)^{1-(\frac{2}{r} + \frac{3}{s})} (\Vert b_{2}\Vert_{L_{t}^{r}L_{x}^{s}} + \Vert b_{3}\Vert_{L_{t}^{r}L_{x}^{s}}) 
     \\&\quad\qquad
     \times
     (\Vert\nabla_{h} u\Vert_{L_{t}^{\infty}L_{x}^2} + \Vert\nabla_{h} b\Vert_{L_{t}^{\infty}L_{x}^2})^{1-\frac{3}{s}} (\Vert\nabla\nabla_{h}u\Vert_{L_{t}^2L_{x}^2} + \Vert\nabla\nabla_{h}b\Vert_{L_{t}^2L_{x}^2})^{1+\frac{3}{s}}
     \\&\qquad
     +
     C(T_2 - T_1)^{1-(\frac{2}{r} + \frac{3}{s})} (\Vert b_{2}\Vert_{L_{t}^{r}L_{x}^{s}} + \Vert b_{3}\Vert_{L_{t}^{r}L_{x}^{s}})
     \\&\quad\qquad
     \times
     (\Vert\partial_{3} u\Vert_{L_{t}^{2}L_{x}^2} + \Vert\partial_{3} b\Vert_{L_{t}^{2}L_{x}^2})^{\frac{s-2}{4s}} (\Vert\partial_{3} u\Vert_{L_{t}^{\infty}L_{x}^2} + \Vert\partial_{3} b\Vert_{L_{t}^{\infty}L_{x}^2})^{\frac{3s-10}{4s}}
     \\&\quad\qquad
     \times 
     (\Vert\nabla\partial_{3}u\Vert_{L_{t}^2L_{x}^2} + \Vert\nabla\partial_{3}b\Vert_{L_{t}^2L_{x}^2})^{\frac{1}{s}} (\Vert\nabla\nabla_{h}u\Vert_{L_{t}^2L_{x}^2} + \Vert\nabla\nabla_{h}b\Vert_{L_{t}^2L_{x}^2})^{1+\frac{2}{s}}
     \\&
     \leq
     C
     +
     C\epsilon J^2(T_2)
     +
     C\epsilon J^{\frac{4}{3}\frac{3s-6}{4s} + 1 + \frac{2}{s}}
     \\&
     \leq
     C
     +
     C\epsilon J^2(T_2).
\end{align*}
Regarding $\tilde{\uppercase\expandafter{\romannumeral3}}$, 
we proceed similarly as in the estimates for $\tilde{\uppercase\expandafter{\romannumeral2}}$.
Namely, we have
\begin{align*}
     \tilde{\uppercase\expandafter{\romannumeral3}}
     &
     =
     \sum_{j, k = 1}^{3}\sum_{i = 1}^{2}\int_{T_1}^{T_2}\int_{\mathbb{T}^3} \partial_{i}u_{j}\partial_{j}b_{k}\partial_{i}b_{k}\,dx\,d\tau
     \\&
     =
     \sum_{i = 1}^{2} \int_{T_1}^{T_2}\int_{\mathbb{T}^3} \partial_{j}u_{1}\partial_{1}b_{1}\partial_{i}b_{1}\,dx\,d\tau
     +
     \sum_{i = 1}^{2}\sum_{k = 2}^{3}\int_{T_1}^{T_2}\int_{\mathbb{T}^3} \partial_{i}u_{1}\partial_{1}b_{k}\partial_{i}b_{k}\,dx\,d\tau
     \\&\qquad
     +
     \sum_{i = 1}^{2}\sum_{k = 1}^{3}\sum_{j = 2}^{3}\int_{T_1}^{T_2}\int_{\mathbb{T}^3} \partial_{i}u_{j}\partial_{j}b_{k}\partial_{i}b_{k}\,dx\,d\tau
     \\&
     =
     \sum_{i = 1}^{2} \int_{T_1}^{T_2}\int_{\mathbb{T}^3} \partial_{j}u_{1}(-b_{2}\partial_{2} - b_{3}\partial_{3})\partial_{i}b_{1}\,dx\,d\tau
     \\&\qquad
     -
     \sum_{i = 1}^{2}\sum_{k = 2}^{3}\int_{T_1}^{T_2}\int_{\mathbb{T}^3} b_{k}\partial_{i}u_{1}\partial_{1i}^{2}b_{k}\,dx\,d\tau
     -
     \sum_{i = 1}^{2}\sum_{k = 2}^{3}\int_{T_1}^{T_2}\int_{\mathbb{T}^3} b_{k}\partial_{1}b_{k}\partial_{ii}^{2}u_{1}\,dx\,d\tau
     \\&\qquad
     -
     \sum_{i = 1}^{2}\sum_{k = 1}^{3}\sum_{j = 2}^{3}\int_{T_1}^{T_2}\int_{\mathbb{T}^3} u_{j}\partial_{j}b_{k}\partial_{ii}^{2}b_{k}\,dx\,d\tau
     -
     \sum_{i = 1}^{2}\sum_{k = 1}^{3}\sum_{j = 2}^{3}\int_{T_1}^{T_2}\int_{\mathbb{T}^3} u_{j}\partial_{i}b_{k}\partial_{ij}^{2}b_{k}\,dx\,d\tau
     \\&
     \leq
     C\int_{T_1}^{T_2}\int_{\mathbb{T}^3} (|u_{2}| + |u_{3}| + |b_{2}| + |b_{3}|) (|\nabla_{h}u| + |\nabla_{h}b| + |\partial_{3}u| + |\partial_{3}b|)(|\nabla\nabla_{h}u| + |\nabla\nabla_{h}b|)\,dx\,d\tau.
\end{align*}
Whence, by H\"older's inequality and Gagliardo-Nirenberg-Sobolev inequality 
the far right side of the above inequality is also bounded by
$$ C + C\epsilon J^2(T_2) + C\epsilon J^{\frac{4}{3}\frac{3s-6}{4s} + 1 + \frac{2}{s}}$$
hence by $C + C\epsilon J^2(T_2)$ in view of (\ref{Ineq01}).
The term $\tilde{\uppercase\expandafter{\romannumeral4}}$ 
is bounded similarly as $\tilde{\uppercase\expandafter{\romannumeral3}}$ 
by $C + C\epsilon J^2(T_2)$, thus, we omit the details.
Next we estimate $\tilde{\uppercase\expandafter{\romannumeral5}}$.
Observing Theorem~\ref{T2}, we have
\begin{align*}
     \tilde{\uppercase\expandafter{\romannumeral5}}
     &
     =
     g\sum_{k= 1}^{3}\sum_{i = 1}^{2}\int_{T_1}^{T_2}\int_{\mathbb{T}^3} \theta e_3\partial_{ii} u_{k}\,dx\,d\tau
     \leq
     C\Vert\theta\Vert_{L_{x, t}^{2}} \Vert\nabla\nabla_{h}u\Vert_{L_{x, t}^{2}}
     \leq
     C\epsilon J(T_2),
\end{align*}
due to (\ref{Ineq}).
Combining the above estimates for $\tilde{\uppercase\expandafter{\romannumeral1}}$ through $\tilde{\uppercase\expandafter{\romannumeral5}}$, we get
\begin{align*}
     &\quad
     \frac{1}{2}\left(\Vert \nabla_{h}u(T_2)\Vert_{L_{x}^2}^2
     +
     \Vert \nabla_{h}b(T_2)\Vert_{L_{x}^2}^2
     \right)
     +
     \int_{T_1}^{T_2} \int_{\mathbb{T}^3} \Vert \nabla\nabla_{h}u\Vert_{L_{x}^2}^2
     +
     \Vert \nabla\nabla_{h}b\Vert_{L_{x}^2}^2\,dx\,d\tau
     \\&
     \leq
     \frac{1}{2}\left(\Vert \nabla_{h}u(T_1)\Vert_{L_{x}^2}^2
     +
     \Vert \nabla_{h}b(T_1)\Vert_{L_{x}^2}^2
     \right)
     +
     C
     +
     C\epsilon J(T_2)
     +
     C\epsilon J^{2}(T_2),
\end{align*}
where is the constant $C$ is independent of $T_2$.
Therefore, we get
\begin{align*}
     \frac{1}{2}J^2(T_2)
     &=
     \sup_{\tau \in (t_1, t_2)} \left\{ \Vert \nabla_{h} u(\tau) \Vert_{2}^2
     +
     \Vert \nabla_{h} b(\tau) \Vert_{2}^2
     \right\}
     +
     \int_{t_1}^{t_2} \Vert \nabla \nabla_{h} u(\tau)\Vert_{2}^2
     +
     \Vert \nabla \nabla_{h} b(\tau)\Vert_{2}^2\,dx\,d\tau
     \\&
     \leq
     \frac{1}{2}\left(\Vert \nabla_{h}u(T_1)\Vert_{L_{x}^2}^2
     +
     \Vert \nabla_{h}b(T_1)\Vert_{L_{x}^2}^2
     \right)
     +
     C\epsilon J(T_2)
     +
     C\epsilon J^{2}(T_2)
     +
     C,
\end{align*}
where we applied the $\epsilon$-Young inequality.
Hence, by choosing $\epsilon < 1/4C$ we obtain
\begin{align}
     &\quad
     \frac{1}{4}\sup_{\tau \in (t_1, t_2)} \left\{ \Vert \nabla_{h} u(\tau) \Vert_{2}^2
     +
     \Vert \nabla_{h} b(\tau) \Vert_{2}^2
     \right\}
     +
     \int_{t_1}^{t_2} \Vert \nabla \nabla_{h} u(\tau)\Vert_{2}^2
     +
     \Vert \nabla \nabla_{h} b(\tau)\Vert_{2}^2\,dx\,d\tau
     \\&
     \leq
     \frac{1}{2}\left(\Vert \nabla_{h}u(T_1)\Vert_{L_{x}^2}^2
     +
     \Vert \nabla_{h}b(T_1)\Vert_{L_{x}^2}^2
     \right)
     +
     C,
     \label{Ineq02}
\end{align}
Finally, we have 
\begin{align*}
     \Vert\nabla_{h}u(T_2)\Vert_{L_{x}^2}^2
     +
     \Vert\nabla_{h}b(T_2)\Vert_{L_{x}^2}^2
     \leq
     \frac{1}{2}\left(\Vert \nabla_{h}u(T_1)\Vert_{L_{x}^2}^2
     +
     \Vert \nabla_{h}b(T_1)\Vert_{L_{x}^2}^2
     \right)
     +
     C,
\end{align*}
for any $T_2 \in (T_1, T_{\text{max}})$.
Therefore we have $$\sup_{T_2 \in (T_1, T_{\text{max}})}\Vert\nabla_{h}u(T_2)\Vert_{L_{x}^2}^2 \leq C < \infty,$$ 
and by (\ref{Ineq01}) and (\ref{Ineq02}), we obtain
\begin{align*}
     \sup_{T_2 \in (T_1, T_{\text{max}})} \left(J^2(T_2) + L^2(T_2)\right)
     \leq
     C
     <
     \infty,
\end{align*}
which implies $$u, b \in L_{t}^{\infty}([0, T); H^1\cap V) \cap L_{t}^{2}([0, T); H^2\cap V).$$
Thus, by our arguments in previous sections, $u$ and $b$ are smooth up to time $T$. 
In particular, $u$ and $b$ are bounded in $H^3\cap V$. 
Whence, we multiply the equation for $\theta$ in (\ref{S1}) by $-\Delta\theta$, 
integrate by parts over $\mathbb{T}^3$ and obtain
\begin{align*}
    \frac{d}{dt}\Vert\nabla\theta\Vert_{L_{x}^2}^2
    &=
    \sum_{i, j=1}^3 \int_{\mathbb{T}^3}u_{j}\partial_{j}\theta\partial_{ii}\theta\,dx
    %\\&
    \leq
    C\int_{\mathbb{T}^3} |\nabla u| |\nabla\theta|^2\,dx
    \\&
    \leq
    C\Vert\nabla u\Vert_{L_{x}^{\infty}}\Vert\nabla\theta\Vert_{L_{x}^2}^2
    %\\&
    \leq
    C\Vert u\Vert_{H_{x}^3}\Vert\nabla\theta\Vert_{L_{x}^2}^2,
\end{align*}
where we used $\nabla\cdot u=0$ and the Sobolev embedding $H^3 \hookrightarrow L^{\infty}$.
Integrating in time from $T_1$ to $T_2$ and by the fact that $u$ is bounded in $H^3$ independent of $T_2$, 
we have $\theta \in L_{t}^{\infty}([0, T); H^1\cap V)$ due to Gr\"onwall's inequality.
The proof of Theorem~\ref{T0} is thus complete.
\endproof

\appendix
\section{Results regarding the fully inviscid case}
We provide a proof following a similar argument to the one given for the existence and uniqueness for the three-dimensional Euler equations in \cite{K} and \cite{MP}.

{\smallskip\noindent {\em Proof of Theorem~\ref{T1}.}}
The first part of the proof follows similarly to that of Theorem~\ref{T3} and we use the same notation here,
except that we choose the orthogonal projection $P_{N}$ from $H$ to its subspaces $H_{\sigma}$ generated by the functions 
$$\{e^{2\pi ik\cdot x} \mid |k| = \max{k_{i}} \leq N\},$$
for integer $N > 0$ and $k\in\mathbb{Z}^3$. 
For $u^{N}, b^{N} \in H_{\sigma}$, 
and $\theta^{N}$ and $p^{N}$ in the corresponding projected space for scalar funtions, respectively, 
we consider solutions of the following ODE system,
\begin{align*}
 \left\{
    \begin{aligned}
     &
     \frac{d u^{N}}{dt}
     +
     P_{N}B(u^{N}, u^{N})
     +
     \nabla p^{N}
     =
     P_{N}B(b^{N}, b^{N})
     +
     g\theta^{N}e_3,
     \\&
     \frac{d b^{N}}{dt}
     +
     P_{N}B(u^{N}, b^{N})
     =
     P_{N}B(b^{N}, u^{N}),
     \\&
     \frac{d \theta^{N}}{dt}
     +
     P_{N}\mathcal{B}(u^{N}, \theta^{N})
     =
     0,
    \end{aligned}
  \right.
\end{align*}
where we slightly abuse the notation by using $B$ and $\mathcal{B}$ 
to denote the same type of nonlinear terms as were introduced in Section~\ref{sec2}. 
We show that the limit of the sequence of solutions exists and solves of original system (\ref{S3}). 
First, we observe that the above ODE system has solution for any time $T>0$ 
since all terms but the time derivatives are at least locally Lipschitz continuous. 
In particular, by similar arguments as in Section~\ref{sec3}, the solution remains bounded in 
$L_{t}^{\infty}((0, \tilde{T}); H) \cap L_{t}^{\infty}((0, \tilde{T}); H^{m}\cap V)$ 
for some $\tilde{T}$ depending on the $H^3$-norm of the initial data. 
Next, we show that $(u^{N}, b^{N}, \theta^{N})$ is a Cauchy sequence in $L^2$.
For $N' > N$, by subtracting the corresponding equations for $(u^{N}, b^{N}, \theta^{N})$ and $(u^{N'}, b^{N'}, \theta^{N'})$, we obtain
\begin{align*}
 \left\{
    \begin{aligned}
     \frac{d}{dt}(u^{N} - u^{N'})
     &=
     -P_{N}B(u^{N}, u^{N})
     +
     P_{N'}B(u^{N'}, u^{N'})
     +
     P_{N}B(b^{N}, b^{N})
     \\&\qquad
     -
     P_{N'}B(b^{N'}, b^{N'})
     -
     \nabla(p^{N} - p^{N'})
     +
     g(\theta^{N} - \theta^{N'})e_3,
     \\
     \frac{d}{dt}(b^{N} - b^{N'})
     &=
     -P_{N}B(u^{N}, b^{N})
     +
     P_{N'}B(u^{N'}, b^{N'})
     +
     P_{N}B(b^{N}, u^{N})
     -
     P_{N'}B(b^{N'}, u^{N'}),
     \\
     \frac{d}{dt}(\theta^{N} - \theta^{N'})
     &=
     -P_{N}\mathcal{B}(u^{N}, \theta^{N})
     +
     P_{N'}\mathcal{B}(u^{N'}, \theta^{N'}).
    \end{aligned}
  \right.
\end{align*}
Next, we take the inner product of the above equations with 
$(u^{N} - u^{N'})$, $(b^{N} - b^{N'})$, and $(\theta^{N} - \theta^{N'})$. 
Adding all three equations, and using (\ref{embd1}) and (\ref{embd2}) from Lemma~\ref{PL1}, 
we obtain
\begin{align*}
     &\quad
     \frac{1}{2}\frac{d}{dt}
     \left(
     \Vert u^{N} - u^{N'}\Vert_{L_{x}^2}^2 + \Vert b^{N} - b^{N'}\Vert_{L_{x}^2}^2 + \Vert\theta^{N} - \theta^{N'}\Vert_{L_{x}^2}^2
     \right)
     \\&
     =
     g((u^{N} - u^{N'})e_3)(\theta^{N} - \theta^{N'})
     -
     (P_{N}B(u^{N}, u^{N}), u^{N'})
     -
     (P_{N'}B(u^{N'}, u^{N'}), u^{N})
     -
     (P_{N}B(b^{N}, b^{N}), u^{N'})
     \\&\qquad
     -
     (P_{N'}B(b^{N'}, b^{N'}), u^{N})
     +
     (P_{N}B(u^{N}, b^{N}), b^{N'})
     +
     (P_{N'}B(u^{N'}, b^{N'}), b^{N})
     -
     (P_{N}B(b^{N}, u^{N}), b^{N'})
     \\&\qquad
     -
     (P_{N'}B(b^{N'}, u^{N'}), b^{N})
     +
     (P_{N}\mathcal{B}(u^{N}, \theta^{N}), \theta^{N'})
     -
     (P_{N'}\mathcal{B}(u^{N'}, \theta^{N'}), \theta^{N})
     \\&
     =
     g((u^{N} - u^{N'})e_3)(\theta^{N} - \theta^{N'})
     +
     ((1 - P_{N})B(u^{N}, u^{N}), u^{N'})
     +
     (B(u^{N} - u^{N'}, u^{N'} - u^{N}), u^{N})
     \\&\qquad
     +
     ((1 - P_{N})B(b^{N}, b^{N}), u^{N'})
     +
     (B(b^{N} - b^{N'}, u^{N'} - u^{N}), u^{N})
     \\&\qquad
     +
     ((1 - P_{N})B(b^{N}, u^{N}), b^{N'})
     +
     (B(b^{N} - b^{N'}, b^{N'} - b^{N}), u^{N})
     \\&\qquad
     -
     ((1 - P_{N})B(u^{N}, b^{N}), u^{N'})
     +
     (B(u^{N} - u^{N'}, b^{N'} - b^{N}), b^{N})
     \\&\qquad
     -
     ((1 - P_{N})\mathcal{B}(u^{N}, \theta^{N}), \theta^{N'})
     +
     (\mathcal{B}(u^{N} - u^{N'}, \theta^{N'} - \theta^{N}), \theta^{N})
     \\&
     =
     S + \sum_{i=1}^{10}S_{i},
\end{align*}
where we integrated by parts 
and used the divergence free condition 
$\nabla\cdot u^{N} = \nabla\cdot u^{N'} = \nabla\cdot b^{N} = \nabla\cdot b^{N'} = 0$.
Then we estimate $S$ and the two types of terms $S_{i}$, $i=1, \ldots, 10$ separately. 
After integration by parts, we first have
\begin{align*}
     S + \sum_{i \text{ even}}S_{i}
     &
     \leq
     g\Vert u^{N} - u^{N'}\Vert_{L_{x}^2} \Vert \theta^{N} - \theta^{N'}\Vert_{L_{x}^2}
     +
     \Vert\nabla u^{N}\Vert_{L_{x}^{\infty}} \Vert u^{N} - u^{N'}\Vert_{L_{x}^2}^2
     \\&\qquad
     +
     2\Vert\nabla b^{N}\Vert_{L_{x}^{\infty}} \Vert u^{N} - u^{N'}\Vert_{L_{x}^2} \Vert b^{N} - b^{N'}\Vert_{L_{x}^2}
     +
     \Vert\nabla u^{N}\Vert_{L_{x}^{\infty}} \Vert b^{N} - b^{N'}\Vert_{L_{x}^2}^2
     \\&\qquad
     +
     \Vert\nabla \theta^{N}\Vert_{L_{x}^{\infty}} \Vert u^{N} - u^{N'}\Vert_{L_{x}^2} \Vert \theta^{N} - \theta^{N'}\Vert_{L_{x}^2}
     \\&
     \leq
     C\left(
     \Vert u^{N} - u^{N'}\Vert_{L_{x}^2}^2 + \Vert b^{N} - b^{N'}\Vert_{L_{x}^2}^2 + \Vert \theta^{N} - \theta^{N'}\Vert_{L_{x}^2}^2
     \right),
\end{align*}
where we used H\"older's inequality 
and the Sobolev embedding $H^{3} \hookrightarrow L^{\infty}$.
Here the constant $C$ depends only on the $H^3$ norm of $u_0$, $b_0$, and $\theta_0$.
Regarding the remaining terms, 
we denote by $\hat{f}$, the Fourier transform of $f \in L^2(\mathbb{T}^3)$  
\begin{align*}
\hat{f}(k) = \frac{1}{(2\pi)^{3/2}} \int_{\mathbb{T}^3}e^{-ik\cdot x}f(x)\,dx,
\end{align*}
and obtain
\begin{align*}
     \sum_{i \text{ odd}}S_{i}
     &
     \leq
     \Vert (u^{N}\cdot\nabla)u^{N}\Vert_{L_{x}^2}\Vert (1 - P_{N})u^{N'}\Vert_{L_{x}^2}
     +
     \Vert (b^{N}\cdot\nabla)b^{N}\Vert_{L_{x}^2}\Vert (1 - P_{N})u^{N'}\Vert_{L_{x}^2}
     \\&\qquad
     +
     \Vert (b^{N}\cdot\nabla)u^{N}\Vert_{L_{x}^2}\Vert (1 - P_{N})b^{N'}\Vert_{L_{x}^2}
     +
     \Vert (u^{N}\cdot\nabla)b^{N}\Vert_{L_{x}^2}\Vert (1 - P_{N})b^{N'}\Vert_{L_{x}^2}
     \\&\qquad
     +
     \Vert (u^{N}\cdot\nabla)\theta^{N}\Vert_{L_{x}^2}\Vert (1 - P_{N})\theta^{N'}\Vert_{L_{x}^2}
     \\&
     \leq
     C\Vert\nabla u^{N}\Vert_{L_{x}^{\infty}} \Vert u^{N}\Vert_{L_{x}^2} \left( \sum_{|k|>N}|\hat{u}^{N'}(k)|^2 (1+|k|^2)^3\frac{1}{{(1+N^2)}^3} \right)^{1/2}
     \\&\qquad
     +
     C\Vert\nabla b^{N}\Vert_{L_{x}^{\infty}} \Vert b^{N}\Vert_{L_{x}^2} \left( \sum_{|k|>N}|\hat{u}^{N'}(k)|^2 (1+|k|^2)^3\frac{1}{{(1+N^2)}^3} \right)^{1/2}
     \\&\qquad
     +
     C\Vert\nabla u^{N}\Vert_{L_{x}^{\infty}} \Vert b^{N}\Vert_{L_{x}^2} \left( \sum_{|k|>N}|\hat{b}^{N'}(k)|^2 (1+|k|^2)^3\frac{1}{{(1+N^2)}^3} \right)^{1/2}
     \\&\qquad
     +
     C\Vert\nabla b^{N}\Vert_{L_{x}^{\infty}} \Vert u^{N}\Vert_{L_{x}^2} \left( \sum_{|k|>N}|\hat{b}^{N'}(k)|^2 (1+|k|^2)^3\frac{1}{{(1+N^2)}^3} \right)^{1/2}
     \\&\qquad
     +
     C\Vert\nabla \theta^{N}\Vert_{L_{x}^{\infty}} \Vert u^{N}\Vert_{L_{x}^2} \left( \sum_{|k|>N}|\hat{\theta}^{N'}(k)|^2 (1+|k|^2)^3\frac{1}{{(1+N^2)}^3} \right)^{1/2}
     \\&
     \leq
     \frac{C}{N^3},
\end{align*}
where $C$ depends on the initial datum, 
and we used the fact that 
$$\Vert f\Vert_{H_{x}^{3}} = \sum_{k\in\mathbb{Z}^3} |\hat{f}(k)|^2(1+|k|^2)^3.$$
Summing up the above estimates we have
\begin{align*}
     &\quad
     \frac{d}{dt}
     \left(
     \Vert u^{N} - u^{N'}\Vert_{L_{x}^2}^2 + \Vert b^{N} - b^{N'}\Vert_{L_{x}^2}^2 + \Vert\theta^{N} - \theta^{N'}\Vert_{L_{x}^2}^2
     \right)
     \\&
     \leq
     C\left(
     \Vert u^{N} - u^{N'}\Vert_{L_{x}^2}^2 + \Vert b^{N} - b^{N'}\Vert_{L_{x}^2}^2 + \Vert\theta^{N} - \theta^{N'}\Vert_{L_{x}^2}^2
     \right)
     +
     \frac{C}{N^3},
\end{align*}
which by Gr\"onwall's inequality implies 
$$\Vert u^{N} - u^{N'}\Vert_{L_{x}^2}^2 + \Vert b^{N} - b^{N'}\Vert_{L_{x}^2}^2 + \Vert\theta^{N} - \theta^{N'}\Vert_{L_{x}^2}^2 \leq \frac{C}{N^3}.$$
Sending $N \to \infty$, we obtain the desired Cauchy sequence.
Namely, $(u^{N}, b^{N}, \theta^{N}) \to (u, b, \theta)$ 
with $u, b \in H$ and $\theta \in L_{x}^2$. 
Due to the above convergence and the fact that 
$u^{N}, b^{N} \in H_{x}^3\cap V$ and $\theta \in H_{x}^3$, 
we see that $u$ and $b$ are also bounded in $H_{x}^3\cap V$ 
while $\theta$ is bounded in $H_{x}^3$.
Thus, the existence part of the theorem is proved by easily verifying that $(u, b, \theta)$ 
satisfies system (\ref{S3}) with some pressure $p$ as discussed below. 
In fact, for a test function $\phi(x) \in \mathcal{V}$ and $0 < t < \tilde{T}$, $(u^{N}, b^{N}, \theta^{N})$ satisfies
\begin{align*}
 \left\{
    \begin{aligned}
     &
     (u^{N}(\cdot, t), \phi)
     =
     (u^{N}((\cdot, 0), \phi)
     +
     \int_{0}^{t}(P_{N}((u^{N}\cdot\nabla)\phi, u^{N})\,d\tau
     -
     \int_{0}^{t}(P_{N}((b^{N}\cdot\nabla)\phi), b^{N})\,d\tau
     +
     g\int_{0}^{t}(\theta^{N}e_3, \phi)\,d\tau,
     \\&
     (b^{N}((\cdot, t), \phi)
     =
     (b^{N}((\cdot, 0), \phi)
     +
     \int_{0}^{t}(P_{N}((u^{N}\cdot\nabla)\phi), b^{N})\,d\tau
     -
     \int_{0}^{t}(P_{N}((b^{N}\cdot\nabla)\phi), u^{N})\,d\tau,
     \\&
     (\theta^{N}((\cdot, t), \phi)
     =
     (\theta^{N}((\cdot, 0), \phi)
     +
     \int_{0}^{t}(\mathcal{B}(u^{N}, \phi), \theta^{N}).
    \end{aligned}
  \right.
\end{align*}
Sending $N \to \infty$ and extracting a subsequence if necessary, 
we have that the integrals of nonlinear terms converge weakly 
to the corresponding integrals of nonlinear terms in (\ref{S3}). 
Also, we see that the nonlinear terms are weakly continuous in time.
Whence by differentiating the first equation in time,  
we conclude that the limit indeed satisfies the equations for $u$ in (\ref{S3}) in the weak sense, i.e.,  
\begin{align*}
     \frac{d}{dt}(u((\cdot, t), \phi)
     =
     -((u\cdot\nabla)u, \phi)
     +((b\cdot\nabla)b, \phi)
     +(g\theta e_3, \phi),
\end{align*}
which in turn implies that there exists some $p \in C([0, \tilde{T}]; H^1)$, such that 
\begin{align*}
     \frac{du}{dt}
     +
     (u\cdot\nabla)u
     +
     \nabla p
     =
     (b\cdot\nabla)b
     +
     g\theta e_3.    
\end{align*}

Regarding uniqueness, 
suppose there are two solutions $(u^{(1)}, b^{(1)}, \theta^{(1)})$ 
and $(u^{(2)}, b^{(2)}, \theta^{(2)})$ with the same initial data $(u_0, b_0, \theta_0)$
for (\ref{S3}).
Subtracting the corresponding equations for the two solutions 
and denoting $\tilde{u}$, $\tilde{b}$, and $\tilde{\theta}$ 
for $u^{(1)} - u^{(2)}$, $b^{(1)} - b^{(2)}$, and $\theta^{(1)} - \theta^{(2)}$, respectively, we obtain
\begin{align*}
 \left\{
    \begin{aligned}
     &
     \frac{\partial\tilde{u}}{\partial t}
     +
     (\tilde{u}\cdot\nabla)u^{(1)}
     +
     (u^{(2)}\cdot\nabla)\tilde{u}
     +
     \nabla\tilde{p}
     =
     (\tilde{b}\cdot\nabla)b^{(1)}
     +
     (b^{(2)}\cdot\nabla)\tilde{b}
     +
     g\tilde{\theta}e_{3},
     \\&
     \frac{\partial\tilde{b}}{\partial t}
     +
     (\tilde{u}\cdot\nabla)b^{(1)}
     +
     (u^{(2)}\cdot\nabla)\tilde{b}
     =
     (\tilde{b}\cdot\nabla)u^{(1)}
     +
     (b^{(2)}\cdot\nabla)\tilde{u},
     \\&
     \frac{\partial\tilde{\theta}}{\partial t}
     +
     (\tilde{u}\cdot\nabla)\theta^{(1)}
     +
     (u^{(2)}\cdot\nabla)\tilde{\theta}
     =
     0,
   \end{aligned}
 \right.
\end{align*}
with $\nabla\cdot\tilde{u} = 0 = \nabla\tilde{b}$ and $\tilde{u}(0)=\tilde{b}(0)=\tilde{\theta}(0)=0$.
Multiply the above equations by $\tilde{u}$, $\tilde{b}$, and $\tilde{\theta}$, respectively, 
integrate over $\mathbb{T}^3$, and add, we get
\begin{align*}
     &\quad
     \frac{1}{2}\frac{d}{dt}
     \left(
     \Vert\tilde{u}\Vert_{L_{x}^2}^2 + \Vert\tilde{b}\Vert_{L_{x}^2}^2 + \Vert\tilde{\theta}\Vert_{L_{x}^2}^2
     \right)
     \\&
     =
     \int_{\mathbb{T}^3}(\tilde{u}\cdot\nabla)u^{(1)}\tilde{u}\,dx
     -
     \int_{\mathbb{T}^3}(\tilde{b}\cdot\nabla)b^{(1)}\tilde{u}\,dx
     +
     \int_{\mathbb{T}^3}g\tilde{\theta}e_{3}\tilde{u}\,dx
     \\&\qquad
     +
     \int_{\mathbb{T}^3}(\tilde{u}\cdot\nabla)b^{(1)}\tilde{b}\,dx
     -
     \int_{\mathbb{T}^3}(\tilde{b}\cdot\nabla)u^{(1)}\tilde{b}\,dx
     +
     \int_{\mathbb{T}^3}(\tilde{u}\cdot\nabla)\theta^{(1)}\tilde{\theta}\,dx
     \\&
     \leq
     C\Vert u^{(1)}\Vert_{L_{x}^{\infty}}\Vert\tilde{u}\Vert_{L_{x}^2}^2
     +
     C\Vert b^{(1)}\Vert_{L_{x}^{\infty}}\Vert\tilde{u}\Vert_{L_{x}^2}\Vert\tilde{b}\Vert_{L_{x}^2}
     +
     C\Vert u^{(1)}\Vert_{L_{x}^{\infty}}\Vert\tilde{u}\Vert_{L_{x}^2}\Vert\tilde{b}\Vert_{L_{x}^2}
     \\&\qquad
     +
     C\Vert \theta^{(1)}\Vert_{L_{x}^{\infty}}\Vert\tilde{u}\Vert_{L_{x}^2}\Vert\tilde{\theta}\Vert_{L_{x}^2},
\end{align*}
where we applied H\"older's inequality and the Sobolev-Nirenberg inequality.
Now due to the embedding $H^3 \hookrightarrow L^{\infty}(\mathbb{T}^3)$, 
and Young's inequality, we have
\begin{align*}
     &
     \frac{1}{2}\frac{d}{dt}
     \left(
     \Vert\tilde{u}\Vert_{L_{x}^2}^2 + \Vert\tilde{b}\Vert_{L_{x}^2}^2 + \Vert\tilde{\theta}\Vert_{L_{x}^2}^2
     \right)
     \leq
     C\left(
     \Vert\tilde{u}\Vert_{L_{x}^2}^2 + \Vert\tilde{b}\Vert_{L_{x}^2}^2 + \Vert\tilde{\theta}\Vert_{L_{x}^2}^2
     \right),
\end{align*}
where $C$ depends on $g$ and $H^{3}$ norm of $(u^{(1)}, b^{(1)}, \theta^{(1)})$.
Thus, by Gr\"onwall's inequality, $(\tilde{u}(t), \tilde{b}(t), \tilde{\theta}(t))$ remains $0$ for $0 \leq t \leq \bar{T}$.
Uniqueness is proved.

\endproof

%\nnewpage
%\section*{Acknowledgments} The author would like to thank...
\bibliographystyle{elsarticle-num}
\bibliography{MHD_Boussinesq_triple}
%\nnewpage
%\bibliography{MHD_Boussinesq_triple}{}
%\bibliographystyle{plain}

\end{document}